\documentclass[11pt]{amsart}
\usepackage[utf8]{inputenc}
\usepackage[T1]{fontenc}
\usepackage[english]{babel}
\usepackage{amsmath}
\usepackage{amsthm}
\usepackage{amssymb}
\usepackage{mathtools}
\usepackage{stmaryrd}

\usepackage{mlmodern}
\usepackage[left=3cm,right=3cm,top=3cm,bottom=3cm]{geometry}

\usepackage{eucal}

\makeatletter
\def\l@subsection{\@tocline{2}{0pt}{2.5pc}{5pc}{}}
\makeatother

\newcommand{\nocontentsline}[3]{}
\newcommand{\tocless}[2]{\bgroup\let\addcontentsline=\nocontentsline#1{#2}\egroup}

\makeatletter
\renewcommand{\l@section}{\@tocline{1}{0pt}{10pt}{1pc}{\bfseries}}
\makeatother

\usepackage[dvipsnames]{xcolor}
\usepackage[all]{xy}
\usepackage{tikz-cd}
\usepackage{array}
\usetikzlibrary{shapes}
\usetikzlibrary[patterns]
\usepackage{caption}
\usepackage{subcaption}
\usepackage[colorlinks=true, citecolor=MidnightBlue, linkcolor=MidnightBlue]{hyperref}
\usepackage{cleveref}

\usepackage{enumitem}
\setitemize{label=$\triangleright$}
\setlist{topsep=8pt, itemsep=4pt }

\makeatletter
\def\thm@space@setup{%
  \thm@preskip=4\topsep \thm@postskip=\thm@preskip
}
\makeatother

 \newtheorem{theo}{Theorem}[section]
 \newtheorem{prop}[theo]{Proposition}
 \newtheorem{defi}[theo]{Definition}
 \newtheorem{lemme}[theo]{Lemma}
 \newtheorem{coro}[theo]{Corollary}

 \newtheorem*{namedtheorem}{\theoremname}
 \newcommand{\theoremname}{Theorem}
 \newenvironment{named}[1]{\renewcommand{\theoremname}{#1}\begin{namedtheorem}}{\end{namedtheorem}}

 \theoremstyle{remark}
 
 \newtheorem{rk}[theo]{Remark}
 \newtheorem{ex}{Example}

 \newenvironment{demo}{\begin{proof}}{\end{proof}}

\newcommand{\floor}{\node[draw,ellipse, minimum width=1cm, minimum height = 0.6 cm]}

\newcommand{\R}{\mathbb{R}}
\newcommand{\Q}{\mathbb{Q}}
\newcommand{\C}{\mathbb{C}}
\newcommand{\N}{\mathbb{N}}
\newcommand{\Z}{\mathbb{Z}}

\newcommand{\CP}{\mathbb{CP}}

\newcommand{\D}{\mathcal{D}}
\newcommand{\calN}{\mathcal{N}}
\newcommand{\G}{\mathcal{G}}
\newcommand{\T}{\mathcal{T}}
\newcommand{\F}{\mathcal{F}}

\renewcommand{\L}{\mathcal{L}}
\renewcommand{\H}{\mathcal{H}}

\renewcommand{\phi}{\varphi}
\renewcommand{\epsilon}{\varepsilon}

\renewcommand{\geq}{\geqslant}
\renewcommand{\leq}{\leqslant}

\newcommand{\<}{\langle}
\renewcommand{\>}{\rangle}
\newcommand{\vide}{\varnothing}

\renewcommand{\div}{\mathrm{div}}
\newcommand{\GL}{\mathrm{GL}}
\newcommand{\som}{\mathrm{sum}}

\newcommand{\dsum}{\displaystyle\sum}
\newcommand{\dprod}{\displaystyle\prod}

\renewcommand{\tilde}{\widetilde}

\newcommand{\LDelta}{\L_\Delta}

\newcommand{\GDelta}{G_\Delta}
\newcommand{\utilde}{{\tilde{u}}}
\newcommand{\codeg}{\mathrm{codeg}}

\newcommand{\mida}{{\lfloor a/2 \rfloor}}
\newcommand{\Mida}{{\left\lfloor \dfrac{a}{2} \right\rfloor}}
\newcommand{\midaD}{{\lfloor a(\Delta)/2 \rfloor}}

\newcommand{\midab}{{\lfloor (a+b)/2 \rfloor}}
\newcommand{\Midab}{{\left\lfloor \dfrac{a+b}{2} \right\rfloor}}

\newcommand{\smax}{s_{\max}}

\newcommand{\stari}{(\star)_i}

\newcommand{\ie}{i.e. }

\title{Universal polynomials for tropical refined invariants in genus 0}
\author{Gurvan Mével}
\address{CNRS \& Nantes Université, UMR 6629 Laboratoire de Mathématiques Jean Leray, 2 rue de la Houssinière, F-44322 Nantes Cedex 3, France}
\email{gurvan.mevel@univ-nantes.fr}

\begin{document}

\begin{abstract}
In \cite{brugalle_polynomiality_2022} the authors showed that the coefficients of small codegree of the tropical refined invariant are polynomial in the Newton polygon. This raised the question of the existence of universal polynomials giving these coefficients, \ie polynomials depending only on the genus and the codegree, and with variables the combinatorial data of the Newton polygon.

In this paper we show that such universal polynomials exist for rational enumeration, and we give an explicit formula. The proof relies on the manipulation of floor diagrams. 
\end{abstract}

\maketitle

\tableofcontents
	
		In this text, by \emph{fan} we mean complete fan in $\R^2$ which is rational with respect to the lattice $\Z^2$. Besides, by \emph{polygon} we mean convex polygon in $\R^2$ with integer vertices.

\section{Introduction}

\subsection{What is it about ?}

\subsubsection{The Göttsche conjecture}

A classical problem in enumerative geometry is to determine the number $N^\delta(d)$ of curves having fixed degree $d$ and $\delta$ nodes, and passing through an appropriate number of points in $\CP^2$. The question can be generalized to smooth surfaces : given a surface $X$, a sufficiently ample line bundle $\L$ over $X$ and $\delta \in \N$, what is the number $N^{X,\delta}(\L)$ of $\delta$-nodal irreducible curves in the linear system $|\L|$ passing through $\frac{\L^2 + c_1(X)\cdot\L}{2}-\delta$ points in generic position ?

For fixed $\delta$, Di Francesco and Itzykson conjectured in \cite{di_francesco_quantum_1995} the number $N^\delta(d)$ to be polynomial for $d$ large enough. For $\delta =1,2,3$ these node polynomials were known in the second half of the XIXth century. They have been computed up to $\delta=14$ by Vainsencher, Kleiman-Piene and Block in \cite{vainsencher_enumeration_1995, kleiman_node_2004, block_computing_2011}.

In \cite{gottsche_conjectural_1998} Göttsche generalized the conjecture of Di Francesco-Itzykson : for all $\delta \in\N$, there exists $P_\delta \in\C[x,y,z,t]$ such that 
for all non-singular complex algebraic surface $X$ and for all line bundle $\L$ sufficiently ample :
\[ N^{X,\delta}(\L) = P_\delta(\L^2, c_1(X)\cdot \L, c_1(X)^2, c_2(X)). \]
Moreover, the generating series of the $P_\delta$'s was conjectured to be multiplicative, \ie there exist some universal power series $A_1,\dots,A_4 \in \Q[[u]] $ such that 
\[ \dsum_{\delta\geq 0} P_\delta(x,y,z,t) u^\delta = A_1^x A_2^y A_3^z A_4^t. \]
This has been proved first by Tzeng \cite{tzeng_proof_2012}. Kool, Shende and Thomas then gave in \cite{kool_short_2011} an alternative proof.

By the adjunction formula one has $g+\delta = \frac{\L^2 - c_1(X)\cdot \L +2}{2}$, where $g$ is the geometric genus of the curve. Hence we can consider a dual problem : given $X,\L,g$, what is $N_{X,g}(\L)$ the number of curves of genus $g$ in $|\L|$ passing through a generic configuration of $c_1(X)\cdot \L -1+g$ points on $X$ ?
However, Di Francesco and Itzykson showed the asymptotic $\ln(N_{\CP^2,0}(d)) \sim 3d \ln(d)$ in \cite{di_francesco_quantum_1995}, hence we cannot hope for the numbers $N_{X,g}(\L)$ to behave polynomially when $g$ is fixed and $\L$ varies.

\subsubsection{Floor diagrams}

Using tropical geometry, Brugallé and Mikhalkin gave in \cite{brugalle_enumeration_2007, brugalle_floor_2008} a combinatorial method to compute these numbers for a certain class of toric surfaces. Via Mikhalkin correspondence theorem \cite{mikhalkin_enumerative_2005} they reduced the enumeration of algebraic curves to the enumeration of a certain type of graphs called floor diagrams.
If the polygon $\Delta$ defines the line bundle $\LDelta$ over the toric surface $X_\F$, where $\F$ is the dual fan of $\Delta$, then the number of curves in $|\LDelta|$ passing through the appropriate number of points can be calculated by counting some marked floor diagrams instead.

 Building on this tool, Fomin and Mikhalkin gave a combinatorial proof of Di Francesco and Itzykson conjecture for $\CP^2$ in \cite{fomin_labeled_2010}.
 Ardila and Block proved polynomiality in \cite{ardila_universal_2013} for families of toric surfaces. Their combinatorial approach allows to also deal with singular surfaces for which the initial Göttsche conjecture does not say anything.
In \cite{block_computing_2014}, Block, Colley and Kennedy considered a logarithmic version of a quantity introduced by Fomin and Mikhalkin and showed it is linear. This gives a new proof for the multiplicativity stated in Göttsche's conjecture in the case of $\CP^2$. Motivated by this work, Liu recovered in \cite{liu_combinatorial_2016} the  result of Block, Colley and Kennedy as a particular case of a more general theorem. Thanks to this theorem, Liu and Osserman completed in \cite{liu_severi_2018} the work of Ardila and Block, especially by clarifying the link with Göttsche conjecture and generalizing the statement to singular surfaces.

\subsubsection{Tropical refined invariants}

Tropical refined invariant for a toric surface $X_\F$ has been introduced by Block and Göttsche in \cite{block_refined_2016}. It is a Laurent polynomial in a variable $q$ which interpolates between complex and real enumeration of curves : plugging $q=1$ we get Gromov-Witten invariant, and plugging $q=-1$ we get tropical Welschinger invariant.
In the case $g=0$ Göttsche and Schroeter defined a tropical refined descendant invariant $\GDelta(s)(q)$ in \cite{gottsche_refined_2019}. Plugging $q=-1$, if $X_\F$ is an unnodal Del Pezzo surface we now get the number of real curves passing through a generic configuration of points having $s$ pairs of complex conjugated points.

In \cite{block_refined_2016} it is shown that if we fix the number of nodes $\delta$, then the coefficients of the tropical refined invariant are eventually polynomial with respect to $\Delta$ : the polynomiality behaviour stated in Göttsche conjecture passes down to the refined level.
More surprisingly, if we fix the genus $g$ instead of the number of nodes $\delta$, then Brugallé and Jaramillo-Puentes proved in \cite{brugalle_polynomiality_2022} that we recover a polynomial behaviour : the coefficients of small codegrees of the tropical refined invariant are eventually given by some polynomials in $\Delta$. This is also the case for the Göttsche-Schroeter invariant : the coefficient of small codegrees of $\GDelta(s)$ are eventually polynomial in $\Delta$ and $s$.

This raises the question of the existence of universal polynomials for these coefficients. 
Is there a Göttsche-like conjecture in this dual and refined setting ? In this paper we positively answer this question for the Göttsche-Schroeter invariant.

\subsection{Notations and results}

	Let $\F$ be a fan and $D(\F)$ the set of all the polygons dual to $\F$. 
	For $i\in\N$ and $\Delta$ a polygon, we 	write $\Delta>i$ if any edge of $\Delta$ has a lattice length greater than $i$. 
	The fan $\F$ determines a toric surface $X_\F$ and any $\Delta \in D(\F)$ gives a line bundle $\LDelta$ and a linear system on $X_\F$. 
	If $\Delta$ is $h$-transverse (see definition \ref{def_polygon}) we will denote by $\GDelta(s)$ the tropical refined invariant (see theorem \ref{def_ITR} which can be taken as a definition).
	If $P$ is a (Laurent) polynomial, its degree $\deg(P)$ is the maximum exponent appearing with a non-zero coefficient, and $\< P\>_i$ is its coefficient of codegree $i$, \ie its coefficient of degree $\deg(P)-i$.

We show in this paper that if $\Delta$ is $h$-transverse, non-singular and large enough with respect to $i$, then $\< \GDelta(s) \>_i$ is given by a polynomial in $s$ and in the combinatorial data of $\Delta$, \ie in the topological data of $(X_\F,\LDelta)$. Moreover, the generating series of these polynomials can be expressed simply and does not depend on the underlying toric surface, making the polynomials \emph{universal}.

To state the result, consider 
\[ A_0 = \dfrac{1}{1-x^2},\ A_1 = \dfrac{1}{1-x},\ A_2 = \dsum_{i\geq0} p(i) x^i = \dprod_{k\geq 1} \dfrac{1}{1-x^k} 
 \]
 where $p(i)$ is the number of partitions of $i$, \ie the number of decreasing sequences $\lambda = (\lambda_1,\dots,\lambda_k)$ of positive integers whose sum is $i$.
  Given three integers $y$, $\chi$ and $s$ we consider the following series
 \[ \G(y,\chi,s)(x) = A_0^s A_1^{y-2-2s} A_2^\chi \]
and we denote by $(P_i)_i$ its coefficients, \ie $\G(y,\chi,s)(x) = \sum_i P_i(y,\chi,s) x^i$. 
For any series $A = \sum_i a_ix^i$ with $a_0=1$ and any integer $n$, the degree $i$ coefficient of $A^n$ is
\[ \dsum_{k_1 + 2k_2 + \dots + ik_i = i}  \binom{n}{k_1,\dots,k_i} a_1^{k_1} \dots a_i^{k_i} . \]
Each term of the sum is polynomial of degree $k_1+\dots+k_i$ in $n$, and thus the degree $i$ coefficient of $A^n$ is polynomial of degree $i$ in $n$.
Hence $P_i$ is polynomial and has degree $i$ in each of the variables $y$, $\chi$ and $s$.
 Given $\Delta$ a polygon with dual fan $\F$ we denote 
 \[ y(\Delta) = c_1(X_\F) \cdot \LDelta \  \text{ and } \  \chi(\Delta) = c_2(X_\F) . \] 
 Combinatorially, these quantities are 
 \[ y(\Delta) = |\partial \Delta \cap \Z^2| \  \text{ and } \  \chi(\Delta) \text{ is the number of vertices of } \Delta . \] 
Our main universality result is the following. We refer to section \ref{subsec_polydiag} for the notations.

\begin{named}{Theorem \ref{theo_nonsing}}
Let $i\in\N$ and $\F$ be a non-singular and $h$-transverse fan. Let $\Delta \in D(\F)$ and $s\in\{0,\dots, \smax(\Delta)\}$. 
If 
\begin{enumerate}
	\item the fan $\F$ has two vertical rays and 
	\[ \left\{ \begin{array}{rcl}
	\Delta &>& 2(i+2) \\
	e^{+\infty}(\Delta) &>& i+ (a(\Delta)-\midaD+1)d_\F \\
	e^{-\infty}(\Delta) &>& \max(i+2s, i+ (a(\Delta)-\midaD+1)d_\F)
	\end{array} \right. \]
\end{enumerate}
or
\begin{enumerate}[resume]
	\item the fan $\F$ has a single vertical ray generated by $(0,\epsilon)$, with $\epsilon \in \{-1,1\}$, and 
	\[ \left\{ \begin{array}{rcl}
	\Delta &>& 2(i+2) \\
	a(\Delta) &>& \max(i+2s, 5(i+1)+6) \\
	e^{\epsilon\infty}(\Delta) &>& \max(i+2s, 5(i+1)+6, i+ (a(\Delta)-\midaD+1)d_\F)
	\end{array} \right. \]
\end{enumerate}
then one has
\[ \< \GDelta(s)\>_i = P_i(y(\Delta),\chi(\Delta),s) .\]
\end{named}

	The strategy of the proof in the following. In the proof of \cite[theorem 1.5]{brugalle_polynomiality_2022} the authors write an explicit formula for $\< G_\Delta(s) \>_i$, where $\Delta$ is a polygon dual to a fan in a certain family. This formula involves floor diagrams of codegree at most $i$ (see section \ref{subsec_codeg}) and constant divergence (see definition \ref{def_diag}).
	 We adapt their proof to a more general family of fans. In particular the divergence of our floor diagrams is no more constant, and we study in lemmas \ref{lemme_r} and \ref{lemme_rgamma} its contribution to the codegree of the floor diagrams. Then we determine the generating series of the different terms appearing in the formula to conclude for this family, and we get theorem \ref{theo_nonsinghori} which corresponds to the point $(1)$ of theorem \ref{theo_nonsing}.
	 
	 For the point $(2)$, we study in proposition \ref{prop_blowup} the link between  $\< \GDelta(s)\>_i$ and $\< G_{\tilde\Delta}(0;s)\>_i$, where the toric surface associated to $\tilde\Delta$ is a toric blow-up of the toric surface associated to $\Delta$. More precisely we construct a correspondence between the floor diagrams with Newton polygon $\Delta$ and the ones with Newton polygon $\tilde\Delta$. 
	 Using this correspondence we are able in corollary \ref{coro_CP2} to include $\CP^2$ in the universality result, but also to lighten the hypothesis of point $(1)$, leading to point $(2)$ of theorem \ref{theo_nonsing}.

With a slight adapation of the proof we can extend the universality result of theorem \ref{theo_nonsinghori} to singular surfaces if we take into account the number $n_k(\Delta)$ of vertices of $\Delta$ of index $k$ (see definition \ref{def_polygon}), \ie the number of singularities of index $k$ of $X_\F$ (see \cite[proposition A.1]{liu_severi_2018}). 
 Given integer $y,s,n_1,n_2,\dots$, we consider the series 
\[ \H(y,s,n_1,n_2,\dots)(x) = A_0^s A_1^{y-2-2s} \dprod_{k\geq 1} A_2(x^k)^{n_k} . \]
We denote by $Q_i$ its degree $i$ coefficient, \ie $\H(y,s,n_1,n_2,\dots)(x) = \sum_i Q_i(y,s,n_1,n_2,\dots) x^i$. Note that $Q_i$ is a polynomial of degree $i$ in the variable $y$ and $s$, and of degree at most $i/k$ in each of the variables $n_k$.

\begin{named}{Theorem \ref{theo_singhori}}
Let $i\in\N$. and $\F$ be a $h$-transverse fan having two vertical rays. Let $\Delta \in D(\F)$ and $s\in \{0,\dots, \smax(\Delta) \}$. If $(\Delta,s)$ satisfies 
\[ \left\{ \begin{array}{rcl}
 \Delta &>& 2(i+2) \\
 e^{-\infty}(\Delta) &>& \max(i+ (a(\Delta)-\midaD+1)d_\F,i+2s) \\
 e^{+\infty}(\Delta) &>& i+ (a(\Delta)-\midaD+1)d_\F
\end{array} \right. \]
 then 
\[ \< \GDelta(s)\>_i = Q_i(y(\Delta),s,n_1(\Delta),\dots). \]
\end{named}

Note that the term $\prod_{k\geq2} A_2(x^k)^{n_k}$, which takes into account the singularities, is the same as the one appearing in \cite[corollary 1.10]{liu_severi_2018}. This is surprising, especially because Liu and Osserman work at fixed number of nodes while we work at fixed genus. It may worth investigate more precisely this phenomenon.

We do not know how to remove the hypothesis on the vertical rays because a singular fan is not as constraint as a non-singular one. As in the proof of theorem \ref{theo_nonsing} we are able to describe what happens when we perform an operation similar to a blow-up, but it is not sufficient to reach any non-singular $h$-transverse fan. 

The paper is organized as follows. In section \ref{sec_diag} we introduce floor diagrams and related quantities. We then prove some lemmas about codegrees that will be used thereafter. Section \ref{sec_univseries} contains the proofs of our universality theorems. In section \ref{sec_removedenom} we briefly discuss what happens if we change the multiplicities of the floor diagrams, leading to simpler formulas. Last, section \ref{sec_lemmas} contains some technical lemmas we used in section \ref{sec_univseries}.

\subsection{Some computations}

The formulas we obtain in theorems \ref{theo_nonsing} and \ref{theo_singhori} are explicit and we can make some computations. 

\begin{ex}
The first coefficients of $\G(y,\chi,s)(x) = A_0^s A_1^{y-2-2s} A_2^\chi$ are
\begin{align*}
	 & P_0(y,\chi,s) = 1, \\
	 & P_1(y,\chi,s) = y + \chi - 2s - 2, \\	 
	 & P_2(y,\chi,s) = \dfrac{1}{2}(y^2 + 2y\chi + \chi^2 - 4ys - 4\chi s + 4s^2 - 3y - \chi + 8s + 2), \\	 
	 & P_3(y,\chi,s) = \dfrac{1}{3!}(y^3 + 3y^2\chi + 3y\chi^2 + \chi^3 - 6y^2 s - 12y\chi s - 6\chi^2 s + 12ys^2 + 12\chi s^2 - 8s^3  \\
	 & \qquad\qquad\qquad\qquad\qquad\qquad\qquad\qquad  - 3y^2 + 3\chi^2 + 18ys + 6\chi s - 24s^2 + 2y - 4\chi - 16s), \\
	  & P_4(y,\chi,s) = \dfrac{1}{4!}(y^4 + 4y^3\chi + 6y^2\chi^2 + 4y\chi^3 + \chi^4 - 8y^3 s - 24y^2\chi s - 24y\chi^2 s - 8\chi^3 s \\
	& \qquad\qquad\qquad\qquad  + 24y^2 s^2 + 48y\chi s^2 + 24\chi^2 s^2 - 32ys^3 - 32\chi s^3 + 16s^4 - 2y^3 + 6y^2\chi \\
	& \qquad\qquad\qquad\qquad + 18y\chi^2 + 10\chi^3 + 24y^2 s - 24\chi^2 s - 72ys^2 - 24\chi s^2 + 64s^3 - y^2 - 14y\chi \\
	& \qquad\qquad\qquad\qquad - \chi^2 - 32ys + 16\chi s + 80s^2 + 2y + 14\chi + 32s).
\end{align*}
\end{ex}

\begin{ex}
In the case of degree $d$ curves on the projective plane $\CP^2$, one has $y(\Delta) = 3d$ and $\chi(\Delta)=3$, see figure \ref{fig_polyCP2}. With corollary \ref{coro_CP2} one has
\begin{figure}[h] 
	\centering
	\begin{tikzpicture}[scale=1] 
		\draw (0,0) node {$\bullet$} ;
		\draw (0,0) node[below] {\scriptsize $(0,0)$} ;
		\draw (2,0) node {$\bullet$} ;
		\draw (2,0) node[below] {\scriptsize $(d,0)$} ;
		\draw (0,2) node {$\bullet$} ; 
		\draw (0,2) node[left] {\scriptsize $(0,d)$} ; 
		\draw (0,0) -- (2,0) ;
		\draw (0,0) -- (0,2) ;
		\draw (2,0) -- (0,2) ;
	\end{tikzpicture}
	\caption{}
	\label{fig_polyCP2}
\end{figure}
\begin{align*}
	& \forall d > \max(11,2s),\ \<G_{\Delta_d}(s)\>_0 = 1, \\	 
	& \forall d > \max(16,2s),\ \<G_{\Delta_d}(s)\>_1 = 3d - 2s + 1, \\
	& \forall d > \max(21,2s),\ \<G_{\Delta_d}(s)\>_2 = \dfrac{1}{2}(9d^2 - 12ds + 4s^2 + 9d - 4s + 8), \\	 
	& \forall d > \max(26,2s),\ \<G_{\Delta_d}(s)\>_3 = \dfrac{1}{3!}(27d^3 - 54d^2 s + 36ds^2 - 8s^3 + 54d^2 - 54ds + 12s^2 \\
	& \qquad\qquad\qquad\qquad\qquad\qquad\qquad\qquad\qquad\qquad\qquad  + 87d - 52s + 42), \\
	& \forall d > \max(31,2s),\ \<G_{\Delta_d}(s)\>_4 = \dfrac{1}{4!}(81d^4 - 216d^3 s + 216d^2 s^2 - 96d s^3 + 16s^4 + 270d^3 \\
	& \qquad\qquad\qquad\qquad\qquad\qquad\qquad\qquad - 432d^2 s + 216ds^2 - 32s^3  + 639d^2 - 744ds + 224s^2 \\
	& \qquad\qquad\qquad\qquad\qquad\qquad\qquad\qquad
	\qquad\qquad\qquad\qquad\qquad\qquad + 690d - 352s + 384).
	 \end{align*}
Note that the bounds we obtain on $d$ in corollary \ref{coro_CP2} are not as sharp as the one in \cite[theorem 1.6 and example 1.9]{brugalle_polynomiality_2022}.
\end{ex}

\tocless{\subsection*{Acknowledgments}}
I am grateful to Erwan Brugallé for his constant support and all his advice regarding this paper and more. 
I would like to thank Assia Mahboubi and Matthieu Piquerez for their willingness to help me when I needed it. I am indebted to Thomas Blomme for his decisive help when it came to handle sums I was unable to manage.  
Part of this paper was written during my stay in Mexico. I thank Cristhian Garay and Lucía López de Medrano for welcoming me, and Benoît Bertrand for a discussion on polygons.

This work was conducted within the France 2030 framework programme, Centre Henri Lebesgue ANR-11-LABX-0020-01. I am supported by the CNRS.

\section{Floor diagrams} \label{sec_diag}

\subsection{$h$-transverse polygons and floor diagrams} \label{subsec_polydiag}

We introduce first some definitions and notations, mainly borrowed from \cite[section 2]{brugalle_polynomiality_2022}.

\begin{defi} \label{def_polygon}
Let $\Delta$ be a polygon and $\F$ be a fan.
\begin{itemize}
\item We said $\Delta$ is \emph{$h$-transverse} if the primitive direction vectors of its edges are of the form $(\pm 1,0)$ or $(n,\pm 1)$ for $n\in\Z$.
Similarly $\F$ is \emph{$h$-transverse} if its rays are generated by vectors of the form $(0,\pm 1)$ or $(\pm 1,n)$ for $n\in\Z$.

\item If $u$ and $v$ are the primitive direction vectors of the edges adjacent to a vertex $P$ of $\Delta$, the \emph{index} of $P$ is $|\det(u,v)|$.  
Similarly, if the $2$-dimensional cone $C$ of $\F$ is generated by rays whose primitive direction vectors are $u$ and $v$, the \emph{index} of $C$ is $|\det(u,v)|$.

\item We said $\Delta$ is \emph{non-singular} if any vertex has index $1$.  
Similarly $\F$ is \emph{non-singular} if any $2$-dimensional cone has index $1$.
\end{itemize} 
\end{defi}

Note that a fan $\F$ is $h$-transverse if any polygon $\Delta \in D(\F)$ is $h$-transverse, or equivalently if there exists a polygon $\Delta \in D(\F)$ which is $h$-transverse. 
Moreover, if $\F$ is a fan and $\Delta \in D(\F)$, the index of a vertex $P$ of $\Delta$ is the index of the dual $2$-dimensional cone $C$ of $\F$. Hence the fan $\F$ is non-singular if any polygon $\Delta \in D(\F)$ is non-singular, or equivalently if there exists a polygon $\Delta \in D(\F)$ which is non-singular. Last, the fan $\F$ has a $2$-dimensional cone $C$ of index $k\geq 2$ if and only if the corresponding toric surface $X_\F$ has a singularity of index $k$ (see \cite[proposition A.1]{liu_severi_2018}). Hence $\F$ is non-singular if and only if the toric surface $X_\F$ is non-singular.

If $\Delta$ is a polygon then we use the notations :
\begin{itemize}
	\item $a(\Delta)$ is the height of $\Delta$, \ie the difference between the maximal and the minimal ordinate of a point of $\Delta$,
	
	\item $e^{+\infty}(\Delta)$ (resp. $e^{-\infty}(\Delta)$) is the length of the top (resp. bottom) horizontal edge of $\Delta$, 	  	
  	\item $y(\Delta) = |\partial \Delta \cap \Z^2|$, 
  	
  	\item $\chi(\Delta)$ is the number of vertices of $\Delta$, \ie the number of $2$-dimensional cones of $\F$, 
  	
  	\item $\smax(\Delta) = \left\lfloor \dfrac{y(\Delta)-1}{2} \right \rfloor$.
\end{itemize}
Note that $y(\Delta) = e^{+\infty}(\Delta) + e^{-\infty}(\Delta) + 2a(\Delta)$.
Moreover, for $\Delta$ a $h$-transverse polygon we denote :
\begin{itemize}
	\item $L(\Delta)$ (resp. $R(\Delta)$) is the unordered list of integers $k\in\Z$ appearing $j$ times, with $j$ maximal such that $j(k,-1)$ is the translation of a edge of the left (resp. right) side of $\Delta$.
\end{itemize}

\begin{ex}
Consider the polygons of figure \ref{fig_expoly}. The polygons $\Delta_1$ and $\Delta_2$ are $h$-transverse but $\Delta_3$ is not, and $\Delta_1$ is non-singular while $\Delta_2$ and $\Delta_3$ are singular. We give in table \ref{table_ex} their combinatorial data.

\begin{figure}[h!]
	\begin{subfigure}[t]{0.33\textwidth}
	\centering
	\begin{tikzpicture}[scale=3/4] 
	
	\foreach \p in {(0,0), (1,0), (2,0), (3,0), (0,1), (1,1), (2,1), (0,2), (1,2), (0,3)}
			{\draw node at \p {$\bullet$} ;}
	\draw (0,0) to (3,0) to (0,3) to cycle;
		
	\end{tikzpicture}
	\caption{$\Delta_1$.}
	\label{fig_expolya}
\end{subfigure}
\begin{subfigure}[t]{0.33\textwidth}
	\centering
	\begin{tikzpicture}[scale=3/4] 
	
	\foreach \p in {(0,0), (1,0), (0,1), (1,1), (2,1), (3,1), (0,2), (1,2), (2,2), (3,2), (1,3), (2,3)}
			{\draw node at \p {$\bullet$} ;}
	\draw (0,0) to (1,0) to (3,1) to (3,2) to (2,3) to (1,3) to (0,2) to cycle;
		
	\end{tikzpicture}
	\caption{$\Delta_2$.}
	\label{fig_expolyb}
\end{subfigure}
\begin{subfigure}[t]{0.32\textwidth}
	\centering
	\begin{tikzpicture}[scale=3/4] 
	
	\foreach \p in {(0,0), (1,0), (2,0), (1,1), (2,1), (3,1), (1,2), (2,2), (3,2), (1,3)}
			{\draw node at \p {$\bullet$} ;}
	\draw (0,0) to (2,0) to (3,1) to (3,2) to (1,3) to cycle;
		
	\end{tikzpicture}
	\caption{$\Delta_3$.}
	\label{fig_expolyc}
\end{subfigure}
	\caption{Some polygons.}
	\label{fig_expoly}
\end{figure}
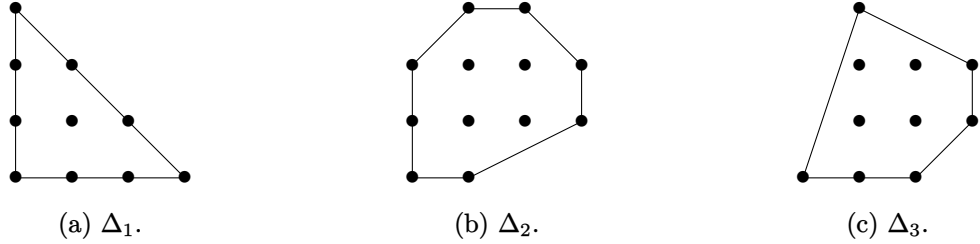

\begin{table}[h!]
\renewcommand{\arraystretch}{1.5}
\centering
\begin{tabular}{|c||c|c|c|c|c|c|c|c|}
\hline 
 & $a(\Delta)$ & $e^{+\infty}(\Delta)$ & $e^{-\infty}(\Delta)$ & $y(\Delta)$ & $\chi(\Delta)$ & $\smax(\Delta)$ & $L(\Delta)$ & $R(\Delta)$ \\
 \hline
 \hline
$\ \Delta_1\ $ & $3$ & $0$ & $3$ & $9$ & $3$ & $4$& $\{ 0,0,0 \}$ & $\{ 1,1,1 \}$  \\
\hline 
$\ \Delta_2\ $ & $3$ & $1$ & $1$ & $8$ & $6$ & $3$ & $\{ -1,0,0 \}$ & $\{ 1,0,-2 \}$ \\
\hline
$\ \Delta_3\ $ & $3$ & $0$ & $2$ & $6$ & $5$ & $2$ & $\not$ & $\not$ \\
\hline 
\end{tabular}
\caption{Combinatorial data of polygons of figure \ref{fig_expoly}.}
\label{table_ex}
\end{table}
\end{ex}

An \textit{oriented graph} $\Gamma$ is a collection of vertices $V(\Gamma)$, of bounded edges $E^0(\Gamma)$ which is a subset of $V(\Gamma)\times V(\Gamma)$ and of infinite edges oriented outward $E^{+\infty}(\Delta)$ (resp. inward $E^{-\infty}(\Gamma)$) which are subsets of $V(\Gamma)$. We denote by $E(\Gamma)$ the set of all edges of $\Gamma$. The graph $\Gamma$ is \textit{weighted} if there is a function $w : E(\Gamma) \to \N^*$. Given a vertex $v\in V(\Gamma)$ of an oriented weighted graph, its \textit{divergence} $\div(v)$ is the difference of the weights entering and leaving $v$. Last, the \emph{genus} of a graph $\Gamma$ is its first Betti number.

\begin{defi}[Floor diagram] \label{def_diag}
Let $\Delta$ be a $h$-transverse polygon. A \emph{floor diagram} $\D$ with Newton polygon $\Delta$ is a quadruple $(\Gamma,w,\ell,r)$ such that 
\begin{itemize}
	\item $(\Gamma,w)$ is a weighted, connected and oriented graph of genus $0$, with $|V(\Gamma)| = a(\Delta)$, $|E^{+\infty}(\Gamma)| = e^{+\infty}(\Delta) $ and $|E^{-\infty}(\Gamma)| = e^{-\infty}(\Delta)$,
	
	\item all the infinite edges have weight 1,
	
	\item $\ell : V(\Gamma) \to L(\Delta) $ and $r : V(\Gamma) \to R(\Delta) $ are bijections such that for every vertex $v \in V(\Gamma)$, $\div(v) = r(v) - \ell(v)$.
\end{itemize}
\end{defi}

\begin{rk}
Floor diagrams are defined in general for any genus. Since we will only consider floor diagrams of genus $0$, we include it in the definition.
\end{rk}

\begin{lemme}
Given a $h$-transverse fan $\F$, there exists $d_\F\in\N$ such that for any floor diagram with Newton polygon $\Delta \in D(\F)$ we have $|\div| \leq d_\F$.
\end{lemme}

\begin{demo}
By hypothesis, the rays of $\F$ are defined by vectors of the form $\pm(0,1)$ or $\pm(1,k)$ for some $k\in\Z$. Let $N_r$ (resp. $n_r$) be the the maximal (resp. minimal) integer $k$ such that $(1,k)$ is a primitive vector of a ray of $\F$. Then the function $r$ of any floor diagram is bounded :
\[ n_r \leq r \leq N_r . \]
Similarly, let $N_\ell$ (resp. $n_\ell$) be the the maximal (resp. minimal) integer $k$ such that $(-1,k)$ is a primitive vector of a ray of $\F$. Then the function $\ell$ of any floor diagram is bounded :
\[ -N_\ell \leq \ell \leq -n_\ell . \]
Since $\div=r-\ell$ one has
\[  n_r + n_\ell \leq \div \leq N_r + N_\ell \]
and the result holds for $d_\F = \max(|N_r+N_\ell|,|n_r+n_\ell|) $.
\end{demo}

By abuse of notations, we will use $\D$ for $\Gamma$. The \textit{degree} of $\D$ is 
\[ \deg(\D) = \dsum_{e\in E(\D)} (w(e)-1) . \]

The orientation of $\D$ induces an partial order $\prec$ on $\D$. We will always draw the floor diagrams oriented from bottom to top. Hence we do not put any arrow on the edges to show the orientation.

\begin{ex}
Figure \ref{fig_exdiag1} gives all the floor diagrams with Newton polygon the polygon of figure \ref{fig_expolya}. The functions $r$ and $\ell$ are constant equal to $1$ and $0$, so any vertex has divergence $1$. 

\begin{figure}[h!]
	\begin{subfigure}[t]{0.33\textwidth}
	\centering
	\begin{tikzpicture}[scale=7/8]
		\floor (3) at (0,4.5) {} ;
		\floor (2) at (0,3) {} ;
		\floor (1) at (-0.3,1.5) {} ;
		
		\draw (3) to (2) ;
		\draw (2) to[out=-115, in=90] (1) ;
		
		\draw (1) to[out=-115, in=90] (-0.6,0);
		\draw (1) to[out=-65, in=90] (0,0);
		\draw (2) to[out=-65, in=90] (0.6,0);

	\end{tikzpicture}
	\caption{}
\end{subfigure}
\begin{subfigure}[t]{0.32\textwidth}
	\centering
	\begin{tikzpicture}[scale=7/8]
		\floor (3) at (0,4.5) {} ;
		\floor (2) at (0,3) {} ;
		\floor (1) at (0,1.5) {} ;
		
		\draw (3) to (2) ;
		\draw (2) to node[left] {\scriptsize $2$} (1) ;		
		
		\draw (1) to[out=-130, in=90] (-0.5,0);
		\draw (1) to[out=-90, in=90] (0,0) ;
		\draw (1) to[out=-50, in=90] (0.5,0);
		
	\end{tikzpicture}
	\caption{}
\end{subfigure}
\begin{subfigure}[t]{0.33\textwidth}
	\centering
	\begin{tikzpicture}[scale=7/8]
		\floor (3) at (1,3) {} ;
		\floor (2) at (-1,3) {} ;
		\floor (1) at (0,1.5) {} ;
		
		\draw (3) to (1) ;
		\draw (2) to (1) ;
		
		\draw (1) to[out=-130, in=90] (-0.5,0);
		\draw (1) to[out=-90, in=90] (0,0) ;
		\draw (1) to[out=-50, in=90] (0.5,0);
	\end{tikzpicture}
	\caption{}
\end{subfigure}
	\caption{The floor diagrams with Newton polygon the polygon of figure \ref{fig_expolya}.}
	\label{fig_exdiag1}
\end{figure}
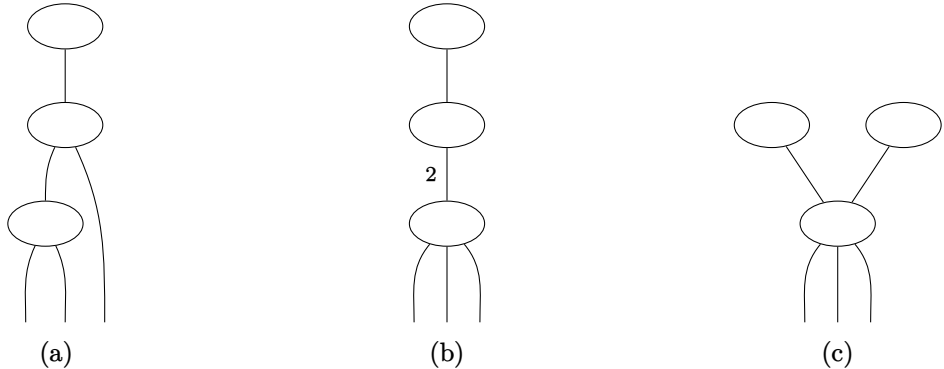
\end{ex}

\begin{ex}
Figure \ref{fig_exdiag2} gives some floor diagrams with Newton polygon the polygon of figure \ref{fig_expolyb}. We show in each vertex the values of $r$ and $\ell$.

\begin{figure}[h!]
	\begin{subfigure}[t]{0.33\textwidth}
	\centering
	\begin{tikzpicture}[scale=7/8]
		
		\floor (3) at (0,4.5) {\scriptsize $-1 \ \ \ \ \ 1$} ;
		\floor (2) at (0,3) {\scriptsize $0 \ \ \ \ \ 0$} ;
		\floor (1) at (0,1.5) {\scriptsize $0 \ \ \ \ \ -2$} ;
		
		\draw (3) to node[left] {\scriptsize $3$} (2) ;
		\draw (2) to node[left] {\scriptsize $3$} (1) ;		
		
		\draw (3) to (0,6) ;
		\draw (1) to (0,0) ;

	\end{tikzpicture}
	\caption{}
\end{subfigure}
\begin{subfigure}[t]{0.32\textwidth}
	\centering
	\begin{tikzpicture}[scale=7/8]
		\floor (3) at (0,4.5) {\scriptsize $0 \ \ \ \ \ 1$} ;
		\floor (2) at (0,3) {\scriptsize $0 \ \ \ \ \ 0$} ;
		\floor (1) at (0,1.5) {\scriptsize $-1 \ \ \ \ \ -2$} ;
		
		\draw (3) to node[left] {\scriptsize $2$} (2) ;	
		\draw (2) to node[left] {\scriptsize $2$} (1) ;	
		
		\draw (1) to (0,0) ;
		
		\draw (3) to (0,6) ;
		
	\end{tikzpicture}
	\caption{}
\end{subfigure}
\begin{subfigure}[t]{0.33\textwidth}
	\centering
	\begin{tikzpicture}[scale=7/8]
		\floor (3) at (1.25,3) {\scriptsize $-1 \ \ \ \ \ 1$} ;
		\floor (2) at (2.5,1.5) {\scriptsize $0 \ \ \ \ \ -2$} ;
		\floor (1) at (0,1.5) {\scriptsize $0 \ \ \ \ \ 0$} ;
		
		\draw (3) to (2) ;	
		\draw (3) to (1) ;	
		
		\draw (1) to (0,0) ;
		
		\draw (2) to (2.5,4.5) ;
		
	\end{tikzpicture}
	\caption{}
\end{subfigure}
	\caption{Some floor diagrams with Newton polygon the polygon of figure \ref{fig_expolyb}.}
	\label{fig_exdiag2}
\end{figure}
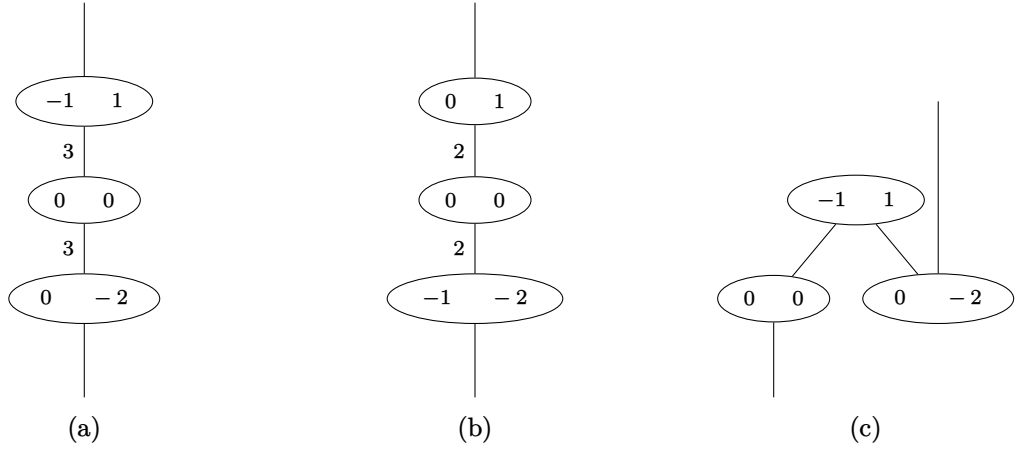
\end{ex}

\begin{defi}[Marking]
Let $\D$ be a floor diagram with Newton polygon $\Delta$. 
A \emph{marking} of $\D$ is an increasing bijection 
\[ m : E(\D) \cup V(\D) \to \{1,\dots,n(\D) \} \]
 where $n(\D)$ is the number of vertices and edges of $\D$. The couple $(\D,m)$ is called a \emph{marked floor diagram}.
 
 Two marked floor diagrams $(\D,m)$ and $(\D',m')$ with Newton polygon $\Delta$ are \emph{isomorphic} if there exists an isomorphism $\phi : \D \to \D'$ of weighted graphs such that $\ell = \ell' \circ \phi$, $r = r'\circ \phi$ and $m=m' \circ \phi$.
 
 We denote by $\nu(\D)$ the number of markings of a diagram $\D$ up to isomorphisms.
\end{defi}

\begin{rk}
A Euler characteristics computation shows that $n(\D) = y(\Delta)-1$.
\end{rk}

\begin{ex}
Figure \ref{fig_exdiagmark} gives examples of marked floor diagrams with Newton polygon the polygon of figure \ref{fig_expolya}. The marked floor diagrams of figures \ref{fig_exdiagmarka} and \ref{fig_exdiagmarkb} are isomorphic.

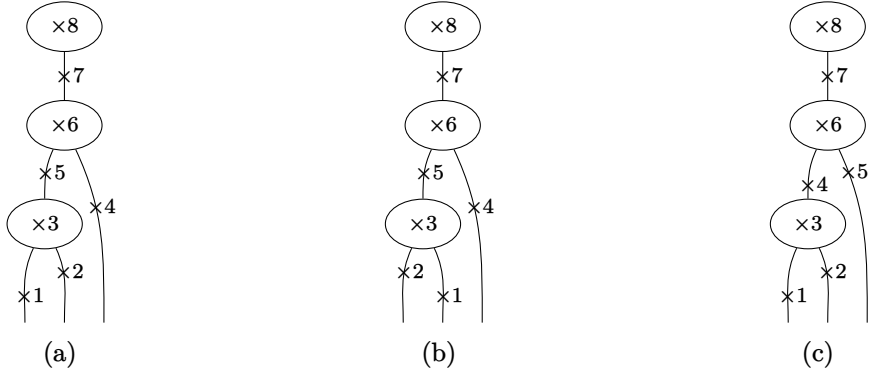
\begin{figure}[h!]
	\begin{subfigure}[t]{0.33\textwidth}
	\centering
	\begin{tikzpicture}[scale=7/8]
		\floor (8) at (0,4.5) {\scriptsize $\times8$};
		\floor (6) at (0,3) {\scriptsize $\times6$};
		\floor (3) at (-0.3,1.5) {\scriptsize $\times3$};
		\draw (8) to node[pos=0.5]{\scriptsize \textcolor{white}{7}$\times7$} (6) ;
		\draw (6) to[out=-115, in=90] node[pos=0.5]{\scriptsize \textcolor{white}{5}$\times5$} (3) ;
		\draw (3) to[out=-115, in=90] node[pos=2/3]{\scriptsize \textcolor{white}{1}$\times1$} (-0.6,0);
		\draw (3) to[out=-65, in=90] node[pos=1/3]{\scriptsize \textcolor{white}{2}$\times2$} (0,0);
		\draw (6) to[out=-65, in=90] node[pos=1/3]{\scriptsize \textcolor{white}{4}$\times4$} (0.6,0);
	\end{tikzpicture}
	\caption{}
	\label{fig_exdiagmarka}
\end{subfigure}
\begin{subfigure}[t]{0.3\textwidth}
	\centering
	\begin{tikzpicture}[scale=7/8]
		\floor (8) at (0,4.5) {\scriptsize $\times8$};
		\floor (6) at (0,3) {\scriptsize $\times6$};
		\floor (3) at (-0.3,1.5) {\scriptsize $\times3$};
		\draw (8) to node[pos=0.5]{\scriptsize \textcolor{white}{7}$\times7$} (6) ;
		\draw (6) to[out=-115, in=90] node[pos=0.5]{\scriptsize \textcolor{white}{5}$\times5$} (3) ;
		\draw (3) to[out=-115, in=90] node[pos=1/3]{\scriptsize \textcolor{white}{2}$\times2$} (-0.6,0);
		\draw (3) to[out=-65, in=90] node[pos=2/3]{\scriptsize \textcolor{white}{1}$\times1$} (0,0);
		\draw (6) to[out=-65, in=90] node[pos=1/3]{\scriptsize \textcolor{white}{4}$\times4$} (0.6,0);
	\end{tikzpicture}
	\caption{}
	\label{fig_exdiagmarkb}
\end{subfigure}	
\begin{subfigure}[t]{0.33\textwidth}
	\centering
	\begin{tikzpicture}[scale=7/8]
		\floor (8) at (0,4.5) {\scriptsize  $\times8$};
		\floor (6) at (0,3) {\scriptsize  $\times6$};
		\floor (3) at (-0.3,1.5) {\scriptsize  $\times3$};
		\draw (8) to node[pos=0.5]{\scriptsize \textcolor{white}{7}$\times7$} (6) ;
		\draw (6) to[out=-115, in=90] node[pos=3/4]{\scriptsize \textcolor{white}{4}$\times4$} (3) ;
		\draw (3) to[out=-115, in=90] node[pos=2/3]{\scriptsize \textcolor{white}{1}$\times1$} (-0.6,0);
		\draw (3) to[out=-65, in=90] node[pos=1/3]{\scriptsize \textcolor{white}{2}$\times2$} (0,0);
		\draw (6) to[out=-65, in=90] node[pos=1/8]{\scriptsize \textcolor{white}{5}$\times5$} (0.6,0);
	\end{tikzpicture}
	\caption{}
	\label{fig_exdiagmarkc}
\end{subfigure}
	\caption{Some marked floor diagrams with Newton polygon the polygon of figure \ref{fig_expolya}.}
	\label{fig_exdiagmark}
\end{figure}
\end{ex}

 A \emph{pairing of order $s$} of the set $P = \{1,\dots,n\}$ is a set $S$ of $s$ disjoint pairs $\{i,i+1\} \subset P$. 
 Given a floor diagram $\D$ and a pairing $S$ of $\{1,\dots,n(\D)\}$, we say that a marking $m$ is \emph{compatible} with $S$ if for any $\alpha \in S$, the set $m^{-1}(\alpha)$ consists of
 \begin{itemize}
 	\item either an edge and an adjacent vertex,
 	\item or two edges that are both entering or both leaving the same vertex.
 \end{itemize}
Let $(\D,m)$ be a marked floor diagram and $S$ a pairing compatible with $m$. We define
\begin{align*}
E_0 &= \{e \in E(\D)\ |\ \forall \alpha\in S, e \notin m^{-1}(\alpha) \}, \\
E_1 &= \{e \in E(\D)\ |\ \exists v\in V(\D), \exists \alpha\in S, \{e,v\} = m^{-1}(\alpha) \} , \\
E_2 &= \{\{e,e'\} \subset E(\D)\ |\ \exists \alpha\in S, \{e,e'\} = m^{-1}(\alpha) \}. \\
\end{align*}
For $n\in\Z$ the quantum integer $[n](q)$ is defined by
\[ [n](q) = \dfrac{q^{n/2}-q^{-n/2}}{q^{1/2}-q^{-1/2}} = q^{(n-1)/2} + q^{(n-3)/2} + \cdots + q^{-(n-3)/2} + q^{-(n-1)/2} \in \N[q^{\pm 1/2}] \]
and we will use the shorcuts
\[ [n]^2 = [n](q)^2 \text{ and } [n]_2 = [n](q^2).  \]

\begin{defi}[Refined $S$-multiplicity] \label{def_mult}
The refined $S$-multiplicity of the marked floor diagram $(\D,m)$ is
\[ \mu_S(\D,m)(q) = \dprod_{e\in E_0} [w(e)]^2 \dprod_{e\in E_1} [w(e)]_2 \dprod_{\{e,e'\} \in E_2} \dfrac{[w(e)][w(e')][w(e)+w(e')]}{[2]} \in \N[q^{\pm 1/2}] \]
if $S$ and $m$ are compatible, and $\mu_S(\D,m)(q)=0$ otherwise.
\end{defi}

\begin{theo}[{\cite[theorem 2.13]{brugalle_polynomiality_2022}}]\label{def_ITR}
Let $\Delta$ be a $h$-transverse polygon and $s\in \N$. For any pairing $S$ of order $s$ of $\{1,\dots, |\mathring\Delta \cap \Z^2|-1 \}$ one has
\[ \GDelta(s) = \dsum_{(\D,m)} \mu_S(\D,m)  \]
where the sum runs over the isomorphism classes of marked floor diagrams with Newton polygon $\Delta$.
\end{theo}

\begin{rk}
The theorem implies that the right-hand side does not depend on the pairing $S$ but only on its order $s$. Thus, to study $\GDelta(s)$ we can choose a particular pairing which makes the calculations easier.
\end{rk}

\begin{rk}
If $\mu_S(\D,m) \neq 0$ then $\deg(\mu_S(\D,m)) = \deg(\D)$.  Let $\D$ be the diagram of figure \ref{fig_degG}, where the function $r$ is increasing and the function $\ell$ is decreasing. The diagram $\D$ has a single marking $m$ and any pairing $S$ is compatible with $m$. We denote by $\mathring\Delta$ the interior of $\Delta$. The weight of the edge $e_k$ is 
\[ w(e_k) = e^{+\infty}(\Delta) + \dsum_{n=k+1}^{a(\Delta)} (r(v_n)-\ell(v_n)) = |\mathring\Delta \cap \{ j=k \} \cap \Z^2| +1  \]
and so $\deg(\D)= |\mathring\Delta \cap \Z^2|$. Because this is the maximal possible degree of a diagram with Newton polygon $\Delta$, we conclude that $\deg(\GDelta(s)) = |\mathring\Delta \cap \Z^2|$.  
\begin{figure}[h!]
	\centering
	\begin{tikzpicture}[scale=7/8]
		
		\floor (4) at (0,5.5) {\scriptsize $v_{a(\Delta)}$} ;
		\floor (3) at (0,3.5) {\scriptsize $v_{k+1}$} ;
		\floor (2) at (0,2) {\scriptsize $v_k$} ;
		\floor (1) at (0,0) {\scriptsize $v_1$} ;
		
		\draw[dashed] (1) to (2) ;
		\draw (2) to node[left] {\scriptsize $w(e_k)$} (3) ;	
		\draw[dashed] (3) to (4) ;
			
		\draw (1) to[out=-130, in=90] (-0.4,-1) ;
		\draw (1) to[out=-50, in=90]  (0.4,-1) ;
		\node at (0,-0.8) {$\cdots$} ;
		
		\draw (4) to[out=130, in=-90] (-0.4,6.5) ;
		\draw (4) to[out=50, in=-90] (0.4,6.5) ;
		\node at (0,6.3) {$\cdots$} ;

		\foreach \p in {(2,2), (3,2), (4,2), (5,2), (6,2), (7,2), (8,2)}
			{\draw node at \p {$\bullet$} ;}
		\draw (3,0) to (5,0) to (7,1) to (8,2) to (8,3) to (7,4) to (5,4) to (3,3) to (2,2) to (2,1) to cycle ;
		
		\draw[->] (9,0) to (9,5) ;
		\node[above] at (9,5) {$j$} ;
		\node at (9,2) {$-$} ;
		\node[right] at (9,2) {$k$} ; 
		\node at (9,0) {$-$} ;
		\node[right] at (9,0) {$0$} ;

	\end{tikzpicture}
	\caption{A floor diagram $\D$ with its Newton polygon $\Delta$.}
	\label{fig_degG} 
\end{figure}
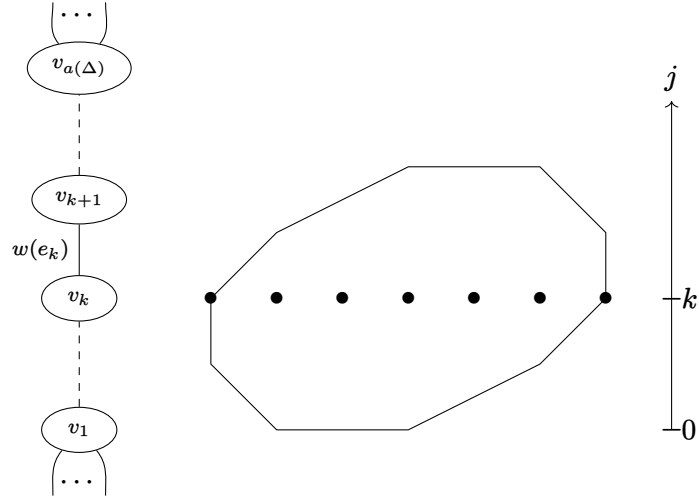
\end{rk}

\begin{rk}
If $\Delta$ and $\Delta'$ are \emph{congruent}, \ie if there exist $A \in \GL_2(\Z)$ and $t$ a translation such that $\Delta' = t(A\Delta)$, then $\GDelta(s) = G_{\Delta'}(s)$.
Indeed, a translation does not change the family of floor diagrams defined by $\Delta$. Moreover, a floor diagram is a way to encode a tropical curve $C$. Via the dual subdivision of $\Delta$ corresponding to $C$, a matrix of $\GL_2(\Z)$ which acts on $\Delta$ also acts on $C$, and preserves its multiplicity. Hence to total count does not change.
\end{rk}

\subsection{Codegrees}\label{subsec_codeg}

Let $\Delta$ be a $h$-transverse polygon and $\D$ a floor diagram with Newton polygon $\Delta$. Its \emph{codegree} is
\[ \codeg(\D) = |\mathring\Delta \cap \Z^2| -\deg(\D) \geq 0. \]
We denote $C_i(\Delta)$ the set of floor diagrams with Newton polygon $\Delta$ and codegree at most $i$.

\begin{rk}
We then have 
\[ \< \GDelta(s) \>_i = \dsum_{(\D,m)} \< \mu_S(\D,m)\>_{i-\codeg(\D)}  \]
where the sum is over the isomorphisms classes of marked floor diagrams with Newton polygon $\Delta$ and codegree at most $i$.
\end{rk}

We will use the following operations on a floor diagram. 
\begin{enumerate}
	\item[$A^+$ : ] If there are $v_1 \prec v_2$ connected by an edge $e_1$ and another edge $e_2$ leaving $v_1$ but not entering $v_2$, then we construct a new diagram as depicted in figure \ref{fig_operationsA+}.
	
	\item[$A^-$ : ] Similarly if $e_2$ is entering $v_2$ but not leaving $v_1$, see figure \ref{fig_operationsA-}.
\end{enumerate}
\begin{figure}[h!]
	\begin{subfigure}[t]{0.49\textwidth}
	\centering
	\begin{tikzpicture}[scale=1]
		\floor (2) at (1,3) {$v_2$} ;
		\floor (1) at (0,1.5) {$v_1$} ;
		
		\draw (1) to node[right,pos=1/2]{\scriptsize $w(e_1)$} (2) ;
		\draw (1) to[out=100, in=-90] node[right,pos=4/5]{\footnotesize $w(e_2)$}  (-0.2,4.5);
		
		\draw[->] (2.2,3) to (3.3,3) ;
		
		\floor (4) at (5,3) {$v_2$} ;
		\floor (3) at (4,1.5) {$v_1$} ;
		
		\draw (3) to node[right,pos=1/2]{\scriptsize $w(e_1)+w(e_2)$} (4) ;
		\draw (4) to[out=90, in=-90] node[left,pos=1/2]{\scriptsize $w(e_2)$}  (5,4.5);
	\end{tikzpicture}
	\caption{Operation $A^+$.}
	\label{fig_operationsA+}
\end{subfigure}
\begin{subfigure}[t]{0.49\textwidth}
	\centering
	\begin{tikzpicture}[scale=1]
		\floor (2) at (0,3) {$v_2$} ;
		\floor (1) at (1,1.5) {$v_1$} ;
		
		\draw (1) to node[right,pos=1/2]{\scriptsize $w(e_1)$} (2) ;
		\draw (2) to[out=-100, in=90] node[right,pos=4/5]{\footnotesize $w(e_2)$}  (-0.2,0);
		
		\draw[->] (2.2,1.5) to (3.3,1.5) ;
		
		\floor (4) at (4,3) {$v_2$} ;
		\floor (3) at (5,1.5) {$v_1$} ;
		
		\draw (3) to node[right,pos=1/2]{\scriptsize $w(e_1)+w(e_2)$} (4) ;
		\draw (3) to[out=-90, in=90] node[left,pos=1/2]{\scriptsize $w(e_2)$}  (5,0);	
	\end{tikzpicture}
	\caption{Operation $A^-$.}
	\label{fig_operationsA-}
\end{subfigure}
	\caption{Operations $A^+$ and $A^-$.}
	\label{fig_operationsA}
\end{figure}
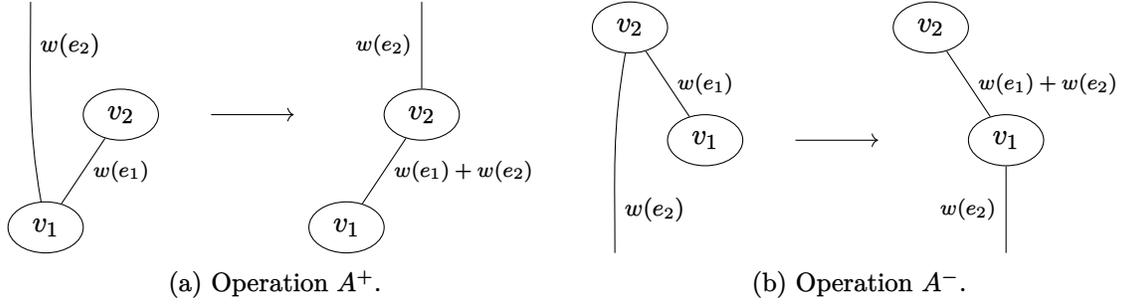
\begin{enumerate}
	\item[$B^\ell$ : ] If there are $v_1 \prec v_2$ connected by an edge $e$ and such that $\ell(v_1) < \ell(v_2)$, then we construct a new diagram as depicted in figure \ref{fig_operationsBl}.
	
	\item[$B^r$ : ] Similarly if $r(v_1) > r(v_2)$, see figure \ref{fig_operationsBr}.
\end{enumerate}
\begin{figure}[h!]
	\begin{subfigure}[t]{0.49\textwidth}
	\centering
	\begin{tikzpicture}[scale=1]
		\floor (2) at (0,3) {\scriptsize $\ell(v_2) \ \ \ $} ;
		\floor (1) at (0,1) {\scriptsize $\ell(v_1) \ \ \ $} ;
		
		\draw (1) to node[right,pos=1/2]{\scriptsize $w(e)$} (2) ;
		
		\draw[->] (1.5,2) to (2,2) ;
		
		\floor (4) at (3,3) {\scriptsize $\ell(v_1) \ \ \ $} ;
		\floor (3) at (3,1) {\scriptsize $\ell(v_2) \ \ \ $} ;
		
		\draw (3) to node[right,pos=1/2]{\scriptsize $w(e) + \ell(v_2) - \ell(v_1)$} (4) ;
	\end{tikzpicture}
	\caption{Operation $B^\ell$.}
	\label{fig_operationsBl}
\end{subfigure}
\begin{subfigure}[t]{0.49\textwidth}
	\centering
	\begin{tikzpicture}[scale=1]
		\floor (2) at (0,3) {\scriptsize $\ \ \ r(v_2)$} ;
		\floor (1) at (0,1) {\scriptsize $\ \ \ r(v_1)$} ;
		
		\draw (1) to node[right,pos=1/2]{\scriptsize $w(e)$} (2) ;;
		
		\draw[->] (1.5,2) to (2,2) ;
		
		\floor (4) at (3,3) {\scriptsize $\ \ \ r(v_1)$} ;
		\floor (3) at (3,1) {\scriptsize $\ \ \ r(v_2)$} ;
		
		\draw (3) to node[right,pos=1/2]{\scriptsize $w(e) + r(v_1) - r(v_2)$} (4) ;
	\end{tikzpicture}
	\caption{Operation $B^r$.}
	\label{fig_operationsBr}
\end{subfigure}
	\caption{Operations $B^\ell$ and $B^r$.}
	\label{fig_operationsB}
\end{figure}
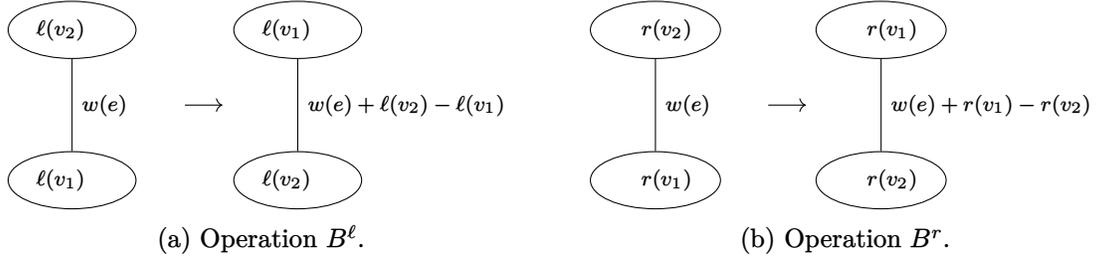
\begin{lemme}[{\cite[lemma 3.2]{brugalle_polynomiality_2022}}]\label{lemme_opeAB}
Genus and Newton polygon are invariant under operations $A^\pm$, $B^{\ell,r}$. Furthermore, the codegree decreases by $w(e_2)$ under operations $A^\pm$, by $\ell(v_2) - l(v_1)$ under operation $B^\ell$ and by $r(v_1) - r(v_2)$ under operation $B^r$.
\end{lemme}

The following lemma is proven in particular cases but with sharper bounds in \cite[lemmas 4.1 and 5.5]{brugalle_polynomiality_2022}.

\begin{lemme}\label{lemmetotalorder}
Let $\Delta$ be a $h$-transverse polygon with dual fan $\F$. 
\begin{enumerate}
 \item[(1)] If $e^{+\infty}(\Delta) > i+d_\F$ then any floor diagram with Newton polygon $\Delta$ and codegree at most $i$ has a unique maximal floor. 
 
 \item[(2)] If $e^{-\infty}(\Delta) > i+d_\F$ then any floor diagram with Newton polygon $\Delta$ and codegree at most $i$ has a unique minimal floor. 
 
 \item[(3)] If $e^{+\infty}(\Delta) > i+d_\F$ and $e^{-\infty}(\Delta) > i+d_\F$ then any floor diagram with Newton polygon $\Delta$ and codegree at most $i$ admits a total order on its set of vertices. 
\end{enumerate} 
\end{lemme}

\begin{demo}
 	The point (3) is an immediat consequence of (1) and (2), 
 	and to prove (2) it suffices to apply (1) to $-\Delta$. Hence we prove (1).

	Let $\D$ be a floor diagram with Newton polygon $\Delta$. Assume that $\D$ admits two maximal vertices. By a finite number of $A^+$ operations we can turn $\D$ into the diagram $\D'$ depicted in figure \ref{fig_2maxvertices_a}.  
	Performing more $A^+$ operations we get the diagram $\D''$ of figure \ref{fig_2maxvertices_b}. By lemma \ref{lemme_opeAB} the codegree decreases by $w_2 + u_1$ under these operations. But $w_2 = e^{+\infty}(\Delta)-u_1+\div(v_2)$ hence the codegree decreases by $e^{+\infty}(\Delta)+\div(v_2)$ and  
\[ \codeg(\D) \geq \codeg(\D'') + e^{+\infty}(\Delta)+\div(v_2) > i. \]
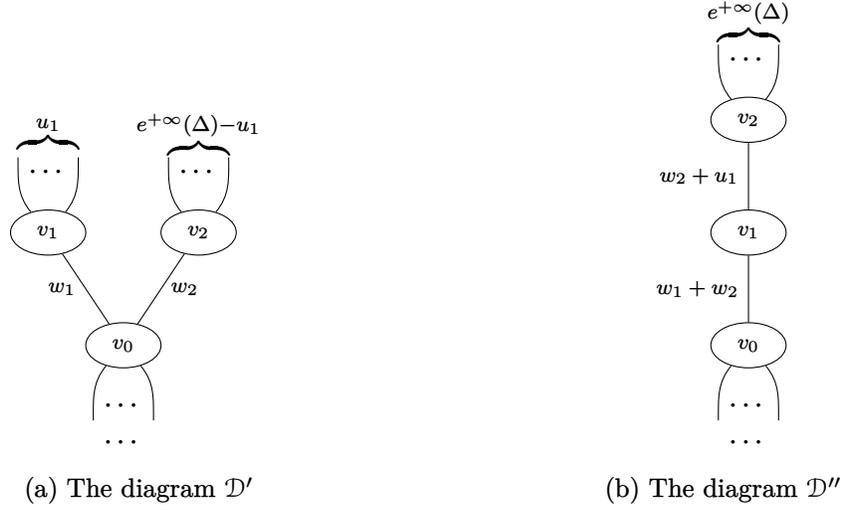
\begin{figure}[h!]
	\begin{subfigure}[t]{0.49\textwidth}
	\centering
	\begin{tikzpicture}[scale=1]
	
	
	\floor (0) at (0,0) {\scriptsize $v_0$} ;
	\floor (1) at (-1,1.5) {\scriptsize $v_1$} ;
	\floor (2) at (1,1.5) {\scriptsize $v_2$} ;
	

	\draw (0) to[out=-130, in=90] (-0.4,-1) ;
	\draw (0) to[out=-50, in=90]  (0.4,-1) ;
	\node at (0,-0.8) {$\cdots$} ;
	\node at (0,-1.3) {$\cdots$} ;
	
	
	\draw (0) to node[left]{\scriptsize $w_1$} (1) ;
	 
	\draw (0) to node[right]{\scriptsize $w_2$} (2) ;
	
	
	\draw (1) to[out=130, in=-90] (-1.4,2.5) ;
	\draw (1) to[out=50, in=-90] (-0.6,2.5) ;
	\node at (-1,2.3) {$\cdots$} ;
	\node at (-1,2.8) {\large $\overbrace{\ \ \ \ \ \ }^{u_1}$} ;
	
	\draw (2) to[out=130, in=-90] (0.6,2.5) ;
	\draw (2) to[out=50, in=-90] (1.4,2.5) ;
	\node at (1,2.3) {$\cdots$} ;
	\node at (1,2.8) {\large $\overbrace{\ \ \ \ \ \ }^{e^{+\infty}(\Delta)-u_1}$} ;

	\end{tikzpicture}
	\caption{The diagram $\D'$}
	\label{fig_2maxvertices_a}
\end{subfigure}	
\begin{subfigure}[t]{0.49\textwidth}
	\centering
	\begin{tikzpicture}[scale=1]	
	
	
	\floor (0) at (0,0) {\scriptsize $v_0$} ;
	\floor (1) at (0,1.5) {\scriptsize $v_1$} ;
	\floor (2) at (0,3) {\scriptsize $v_2$} ;
	

	\draw (0) to[out=-130, in=90] (-0.4,-1) ;
	\draw (0) to[out=-50, in=90]  (0.4,-1) ;
	\node at (0,-0.8) {$\cdots$} ;
	\node at (0,-1.3) {$\cdots$} ;
	
	
	\draw (0) to node[left]{\scriptsize $w_1+w_2$} (1) ;
	 
	\draw (1) to node[left]{\scriptsize $w_2+u_1$} (2) ;
	
		
	\draw (2) to[out=130, in=-90] (-0.4,4) ;
	\draw (2) to[out=50, in=-90] (0.4,4) ;
	\node at (0,3.8) {$\cdots$} ;
	\node at (0,4.3) {\large $\overbrace{\ \ \ \ \ \ }^{e^{+\infty}(\Delta)}$} ;

	\end{tikzpicture}
	\caption{The diagram $\D''$}
	\label{fig_2maxvertices_b}
\end{subfigure}
\caption{Diagram with two maximal vertices.}
	\label{fig_2maxvertices}
\end{figure}
\end{demo}

\begin{lemme}\label{lemme_ordertop}
Let $\D$ be a floor diagram 
 and codegree at most $i$, and $a=|V(\D)|$. If $\D$ has a unique maximal (resp. minimal) floor then the order is total on the $a-i-1$ highest (resp. lowest) floors.
\end{lemme}

\begin{demo}
We will give the proof in the case where $\D$ has a unique top floor, the other point being proved applying that case to $-\Delta$.

Let $b$ be the maximal integer such that the order is total on $v_{a-b} \prec \dots \prec v_a$. We would like to show that $a-b \leq i+2$. Assume the contrary. By $A^\pm$ operations we reduce to the case where $\D$ looks like the diagram of figure \ref{fig_ordertop}. By more operations $A^-$ we can attach $v_c$ under $v_{c+2}$ which reduces the codegree by at least $a-b-c-2$. Then we attach $v_{c+1}$ under $v_1$ by other operations $A^-$ which reduces the codegree by at least $c$. Hence we have
\[ \codeg(\D) > (a-b-c-2) + c = a-b-2 > i  \]
which is a contradiction.
\begin{figure}[h!]
	\centering
	\begin{tikzpicture}[scale=7/8]
	
	
	\floor (0) at (0,2) {$v_a$} ;
	\floor (1) at (0,0) {$v_{a-b}$} ;
	
	\floor (2) at (-1,-1.5) {$v_{a-b-1}$} ;
	\floor (4) at (-1,-3.5) {$v_{c+1}$} ;
	
	\floor (3) at (1,-1.5) {$v_c$} ;
	\floor (5) at (1,-3.5) {$v_1$} ;
	
	
	\draw (0) to[out=130, in=-90] (-0.5,3) ;
	\draw (0) to[out=50, in=-90] (0.5,3) ;
	\node at (0,2.7) {$\cdots$} ;
	
	
	\draw[dashed] (0) to (1) ;
	\draw (1) to (2) ;
	\draw (1) to (3) ;
	\draw[dashed] (2) to (4) ;
	\draw[dashed] (3) to (5) ;

	
	\draw (4) to[out=-130, in=90] (-1.5,-4.5) ;
	\draw (4) to[out=-50, in=90] (-0.5,-4.5) ;
	\node at (-1,-4.2) {$\cdots$} ;

	\draw (5) to[out=-130, in=90] (0.5,-4.5) ;
	\draw (5) to[out=-50, in=90] (1.5,-4.5) ;
	\node at (1,-4.2) {$\cdots$} ;

	\end{tikzpicture}
	\caption{A diagram with a unique maximal floor.}
	\label{fig_ordertop}
\end{figure}
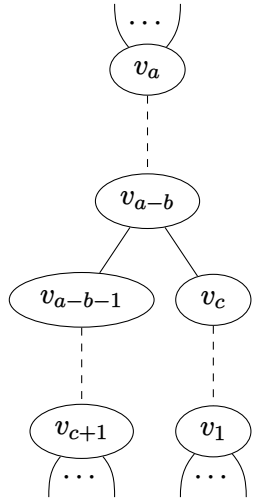
\end{demo}

The following lemma is a generalization of the last assertion of \cite[lemma 4.1]{brugalle_polynomiality_2022}. 

\begin{lemme}\label{lemmepoidsgrand}
	Let $i\in\N$ and $\Delta$ be a $h$-transverse polygon with dual fan $\F$. Let $(a,e^{+\infty},e^{-\infty}) = (a(\Delta),e^{+\infty}(\Delta),e^{-\infty}(\Delta))$ and assume $a>2(i+2)$.
\begin{enumerate}
	\item[(1)] If $e^{+\infty} > i+ (a-\mida+1)d_\F$, then the weights of the $a-\mida$ highest bounded edges of any floor diagram with Newton polygon $\Delta$ and codegree at most $i$ are greater than $i-\codeg(\D)$.
	
	\item[(2)] If $e^{-\infty} > i+ (a-\mida+1)d_\F$, then the weights of the $a-\mida$ lowest bounded edges of any floor diagram with Newton polygon $\Delta$ and codegree at most $i$ are greater than $i-\codeg(\D)$.
	
	\item[(3)] If $e^{\pm\infty} >  i+ (a-\mida+1)d_\F$, then the weights of all the bounded edges of any floor diagram with Newton polygon $\Delta$ and codegree at most $i$ are greater than $i-\codeg(\D)$. 
\end{enumerate}	
	
\end{lemme}

\begin{demo}
	The point (3) is an immediat consequence of (1) and (2), and to prove (2) it suffices to apply (1) to $-\Delta$. Hence we prove (1). 
	Note that the hypotheses imply $e^{+\infty} > i+d_\F$, hence by lemma \ref{lemmetotalorder} any floor diagram $\D$ with Newton polygon $\Delta$ and codegree at most $i$ has a unique maximal floor. By lemma \ref{lemme_ordertop}
 the order is total on the $a-i-1$ highest floors of $\D$, and in particular on its $a-\mida+1$ highest floors. Moreover, if there is an infinite edge attached to $v_k$ for $k< a-i$, then by $A^+$ operations we see that $\codeg(\D) >i$. Thus $\D$ looks like the diagram of figure \ref{fig_totalorder}. 	
\begin{figure}[h!]
	\centering
	\begin{tikzpicture}[scale=7/8]
	
	\floor (2) at (3.5,3) {$v_{i+2}$} ;
	\floor (3) at (3.5,5) {$v_{\mida}$} ;
	\floor (4) at (3.5,7) {$v_{a-i}$} ;
	\floor (5) at (1.5,9) {$v_{a-1}$} ;
	\floor (6) at (0,10.5) {$v_a$} ;
	
	\draw[dashed] (2) to (3) ;
	\draw[dashed] (3) to (4) ;
	\draw[dashed] (4) to (5) ;
	\draw (5) to (6) ;
		
	\draw (2) to (3.5,2) ;
		\node at (3.5,1.5) {$\vdots$} ;

		\draw (6) to[out=120, in=-90] (-0.4,12) ;
		\draw (6) to[out=60, in=-90]  (0.4,12);
		\node at (0,11.4) {$\cdots$} ;
		\node at (0,12.6) {};
		
		\draw (5) to[out=110, in=-90] (1.1,12) ;
		\draw (5) to[out=70, in=-90]  (1.9,12);
		\node at (1.5,11.4) {$\cdots$} ;
		\node at (1.5,12.4) {\large $\overbrace{\ \ \ \ \ }^{u_1}$} ;
		
		\node at (2.5,11.4) {$\cdots$} ;
		
		\draw (4) to[out=105, in=-90] (3.1,12) ;
		\draw (4) to[out=75, in=-90]  (3.9,12);
		\node at (3.5,11.4) {$\cdots$} ;
		\node at (3.5,12.4) {\large $\overbrace{\ \ \ \ \ }^{u_i}$} ;

	\end{tikzpicture}
	\caption{A diagram with a total order on its highest vertices.}
	\label{fig_totalorder}
\end{figure}
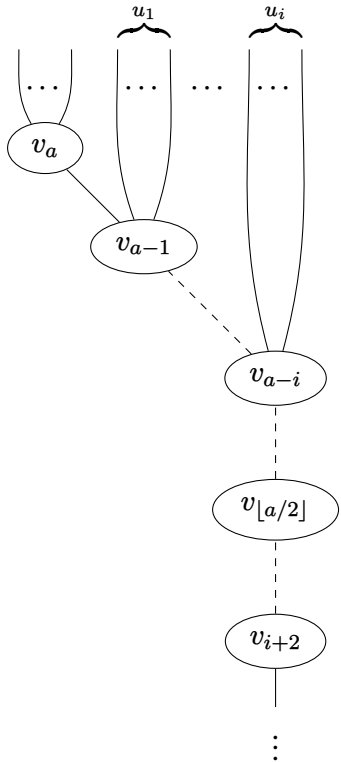		
	For $\mida \leq k \leq a-1$ let $e_k$ be the bounded edge between $v_k$ and $v_{k+1}$.
\begin{itemize} 
	\item If $a-i \leq k \leq a-1$, the weight of $e_k$ is
	\begin{align*}
		w(e_k) &= e^{+\infty} - \dsum_{j=a-k}^i u_j + \dsum_{j=k+1}^a \div(v_j) \\
		&\geq i+\left(a-\Mida+1 \right)d_\F - \codeg(\D) - (a-k)d_\F \\
		&> i-\codeg(\D) + \left( \dfrac{a}{2} +1 -i \right) d_\F \\
		&\geq i-\codeg(\D).
	\end{align*}	 
	In particular $w(e_{a-i}) > i-\codeg(\D) + (a/2+1-i)d_\F$.

	\item If $ \mida \leq k \leq a-i-1$, the weight of $e_k$ is
	\begin{align*}
		w(e_k) &= w(e_{a-i}) + u_i + \dsum_{j=k+1}^{a-i} \div(v_j) \\
		&> i-\codeg(\D) + \left( \dfrac{a}{2}+1-i \right) d_\F - (a-i-k)d_\F \\
		&\geq i-\codeg(\D) + \left( \Mida - \dfrac{a}{2}+1 \right) d_\F \\
		&\geq i-\codeg(\D).
	\end{align*}	 
	
\end{itemize}
\end{demo}

\section{Universal series in genus 0} \label{sec_univseries}

\subsection{The case $e^{\pm\infty}(\Delta) \neq 0$}\label{subsec_maintheopart}

For any integer vector $u\in\N^\N$ (or $\N^{\N^*}$) with finite support and any $k\in\N$ we set 
\[ \som_k(u) = \dsum_{j\geq k} u_j \text{ and } \codeg_k(u) = \dsum_{j\geq k} ju_j.  \]
We will use the shortcut $\codeg = \codeg_1$, and for $i\geq 1$ let
\begin{align*}
C_i &= \{ u \in \N^{\N^*}\ |\ \codeg(u) \leq i\} , \\
B_i &= \{ u \in \N^{\N^*}\ |\ \codeg(u)=i \}. 
\end{align*}
 Note that if $u\in C_i$ or $u\in B_i$, then $u_k=0$ for $k\geq i+1$. Hence we can consider $u$ as a vector in $\N^i$ by forgetting $u_k$ for $k\geq i+1$.
 
For $s\in\N$ we denote by $d(s)$ the set of all decompositions of $s$, \ie
\[ d(s) = \{ S \in \N^\N\ |\ \som_0(S)=s  \}, \]
and for $S \in d(s)$ we set
\[ \binom{s}{S} = \binom{s}{S_0,S_1,\dots} = \dfrac{s!}{S_0! S_1! \dots}. \]

For integers $a,p \in\N$ and vectors $u \in\N^{\N^*}$ and $S\in\N^\N$ we define 
\begin{align*}
\nu_n(a,p,u,S) &= \dbinom{a+np-\som_{n+1}(u-2S)}{u_n-2S_n},\\
\nu_{\geq n}(a,p,u,S) &= \dprod_{k\geq n} \nu_k(a,p,u,S), \\
\calN(a,p,S) &= \dsum_{n\geq0} \left( \dsum_{\codeg(u)=n} \nu_{\geq1}(a,p,u,S)\right) x^n \\
   &= \dsum_{u\in\N^{\N^*}} \nu_{\geq1}(a,p,u,S) x^{\codeg(u)}.
\end{align*} 

Given two integers $k,\ell \geq 0$ we define
\[ F(k,\ell) = \dsum_{\substack{i_1 + \dots + i_k = \ell \\ i_j \geq 1 }} \dprod_{j=1}^k i_j \  \text{ and } \  \Phi_\ell(k) = F(k,k+\ell) \]
with the convention $\Phi_0(0) = 1$.

Recall that we consider the following formal series :
\[
 A_0 = \dfrac{1}{1-x^2},\ A_1 = \dfrac{1}{1-x},\ A_2 = \dsum_{n\geq0} p(n) x^n = \dprod_{k\geq 1} \dfrac{1}{1-x^k} 
 \]
 where $p(n)$ is the number of partitions of $n$.

We postpone to section \ref{sec_lemmas} some lemmas regarding all these quantities that will be used in the proofs of this section.

\begin{lemme} \label{lemme_r}
Let $i\in\N$ and $\Delta$ be a non-singular polygon with $\Delta > 2i$. Let $\D$ be a floor diagram with Newton polygon $\Delta$,  
codegree at most $i$ and having a total order on its vertices $v_1 \prec \dots \prec v_a$. Let $n=\min_k r(v_k)$, $N = \max_k r(v_k)$, and for $n\leq k \leq N$ let $a_k$ be the number of vertices with $r(v)=k$. Finally let $\alpha_k = a_n+\dots + a_k$. One has :
\begin{enumerate}
	\item if $1 \leq j \leq a_n-i$ then $r(v_j) = n$,
	
	\item if $n \leq k \leq N-1$ then :
	\begin{enumerate}
		\item if $\alpha_k -i+1 \leq j \leq \alpha_k+i$ then $r(v_j) \in \{k,k+1\}$ ; moreover there are $i$ vertices with $r(v)=k$ and $i$ vertices with $r(v) = k+1$,
	
	\item if $\alpha_k+i+1 \leq j \leq \alpha_{k+1}-i$ then $r(v_j) = k+1$.
	\end{enumerate}
	
	\item if $\alpha_N-i+1 \leq j \leq \alpha_N$ then $r(v_j) = N$.
\end{enumerate} 
\end{lemme}

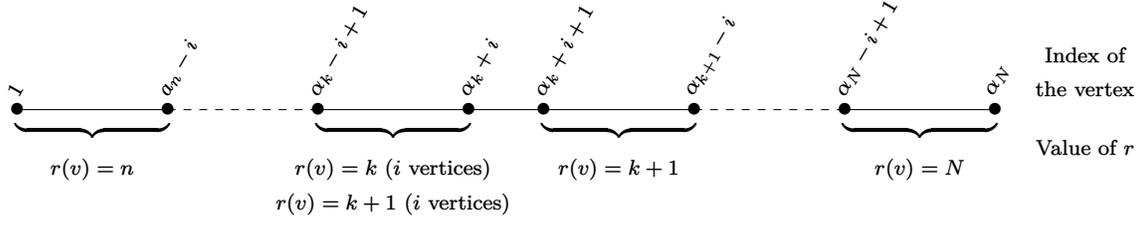
\begin{figure}[h!] 
		\centering
	\begin{tikzpicture}[scale=1]
		
		\draw (-6,0) node {$\bullet$} ;
		\draw (-6,0) node[right, xshift=-5pt,rotate=60] {\scriptsize $\ 1$} ;
		
		\draw (-6,0) to (-4,0) ;	
		
		\draw (-4,0) node {$\bullet$} ;
		\draw (-4,0) node[right, xshift=-5pt,rotate=60] {\scriptsize $\ a_n-i$} ;
		
		\node[below, align=center] at (-5,0) {$\underbrace{\ \ \ \ \ \ \ \ \ \ \ \ \ \ \ \ }_{}$ \\ \scriptsize $r(v) = n$} ;
		
		
		\draw[dashed] (-4,0) to (-2,0) ;
		
		
		\draw (-2,0) node {$\bullet$} ;
		\draw (-2,0) node[right, xshift=-5pt,rotate=60] {\scriptsize $\ \alpha_k-i+1$} ;
		
		\draw (-2,0) to (0,0) ;
		
		\draw (0,0) node {$\bullet$} ;
		\draw (0,0) node[right, xshift=-5pt,rotate=60] {\scriptsize $\ \alpha_k+i$} ;
		
		\node[below, align=center] at (-1,0) {$\underbrace{\ \ \ \ \ \ \ \ \ \ \ \ \ \ \ \ }_{}$ \\ \scriptsize  $r(v) = k$ ($i$ vertices) \\ \scriptsize  $r(v) = k+1$ ($i$ vertices) } ;
		
		
		\draw (0,0) to (1,0) ;
		
		 
		\draw (1,0) node {$\bullet$} ;
		\draw (1,0) node[right, xshift=-5pt,rotate=60] {\scriptsize $\ \alpha_k+i+1$} ;
		
		\draw (1,0) to (3,0) ;
		
		\draw (3,0) node {$\bullet$} ;
		\draw (3,0) node[right, xshift=-5pt,rotate=60] {\scriptsize $\ \alpha_{k+1}-i$} ;
		
		\node[below, align=center] at (2,0) { $\underbrace{\ \ \ \ \ \ \ \ \ \ \ \ \ \ \ \ }_{}$ \\ \scriptsize  $r(v) = k+1$ } ;
		
		
		\draw[dashed] (3,0) to (5,0) ;
		

		\draw (5,0) node {$\bullet$} ;
		\draw (5,0) node[right, xshift=-5pt,rotate=60] {\scriptsize $\ \alpha_N-i+1$} ;
		
		\draw (5,0) to (7,0) ;	
		
		\draw (7,0) node {$\bullet$} ;
		\draw (7,0) node[right, xshift=-5pt,rotate=60] {\scriptsize $\ \alpha_N$} ;
		
		\node[below, align=center] at (6,0) {$\underbrace{\ \ \ \ \ \ \ \ \ \ \ \ \ \ \ \ }_{}$ \\ \scriptsize $r(v) = N$} ;

		\node[align=center] at (8.2,0.5) {\scriptsize  Index of \\ \scriptsize  the vertex} ;	
		\node[align=center] at (8.2,-0.5) {\scriptsize  Value of $r$} ;
	\end{tikzpicture}
	\caption{The function $r$.} 
	\label{fig_gamma}
\end{figure}

The situation described in lemma \ref{lemme_r} is summarized in figure \ref{fig_gamma}.  
 
\begin{demo}
Because $\Delta$ is non-singular, all its vertices have index $1$. Hence the right side of $\Delta$ looks like the picture of figure \ref{fig_rightside}. For any $n\leq k\leq N$, $a_k$ is the length of the edge of $\Delta$ having direction vector $(k,-1)$ so $a_k > 2i$ by hypothesis, which implies $\alpha_k - i > \alpha_{k-1}+i$ and $\alpha_k > 2(k-n+1)i$.
	We investigate how can we choose the function $r$ so that $\D$ has codegree at most $i$. We construct $r$ from the bottom vertices to the top ones. 
	The key element is lemma \ref{lemme_opeAB}.

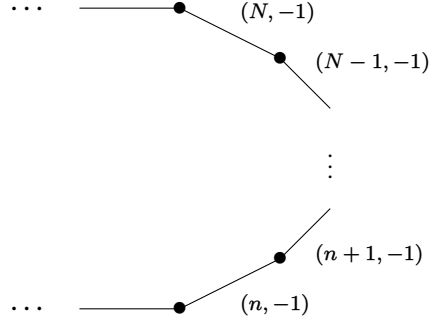
\begin{figure}[h!]
		\centering
	\begin{tikzpicture}[scale=1/1.5]
		
		\node at (-1,0) {$\dots$} ;
		\draw (0,0) to (2,0) to (4,1) to (5,2) ;
		\draw (5,4) to (4,5) to (2,6) to (0,6) ;
		\node at (-1,6) {$\dots$} ;
		
		\node[below right] at (3,0.5) {\scriptsize  $(n,-1)$} ;
		\node[below right] at (4.5,1.5) {\scriptsize $(n+1,-1)$} ;
		\node at (5,3) {\scriptsize $\vdots$} ;
		\node[above right] at (4.5,4.5) {\scriptsize $(N-1,-1)$} ;
		\node[above right] at (3,5.5) {\scriptsize $(N,-1)$} ;

		\node at (2,0) {$\bullet$} ;
		\node at (4,1) {$\bullet$} ;
		
		\node at (4,5) {$\bullet$} ;
		\node at (2,6) {$\bullet$} ;

	\end{tikzpicture}
	\caption{The right side of $\Delta$.}
	\label{fig_rightside}
\end{figure}

\begin{enumerate}
	\item Assume that $r(v_j) > n$ for some $1 \leq j \leq a_n-i$. The vertex $v_j$ has at least $a_n-j+1$ vertices with $r(v)=n$ above it, thus we can perform at least $a_n-j+1$ operations $B^r$, each of them making the codegree drop by at most $1$ by lemma \ref{lemme_opeAB}. Since $-j \geq i-a_n$ we get $\codeg(\D) \geq i+1$, a contradiction.

	\item We prove it by induction over $k$.
	\begin{itemize}
	\item If $k=n$ then necessarily $r(v_j) \geq n$ for any $j$.
	Among the $a_n$ vertices having $r(v)=n$ we know by (1) that $a_n-i$ of them are $v_1,\dots, v_{a_n-i}$. Thus it remains $i$ vertices to be given $r(v)=n$.
	If $r(v_j)=n$ with $j > a_n+i$ then between $v_{a_n-i+1}$ and $v_{a_n+i}$ there are at most $i-1$ vertices having $r(v)=n$, so at least $i+1$ vertices having $r(v) > n$. These $i+1$ vertices are all below $v_j$ thus we can perfom at least $i+1$ operations $B^r$ each of them making the codegree drop by at least $1$. Hence $\codeg(\D) > i$, contradicton.
	
	Thus between $v_{a_n-i+1}$ and $v_{a_n+i}$ there are $i$ vertices with $r(v) = n$. Assume that $r(v_j)>n+1$ for some $a_n-i+1 \leq j \leq a_n+i$. Then $v_j$ has at most $j-a_n+i$ vertices with $r(v) \in \{n, n+1\}$ below it, so at least 
	\[ (i+a_{n+1}) - (j-a_n+i) = a_n+a_{n+1}-j \geq a_{n+1}-i > i\]
	vertices with $r(v) \in \{n, n+1\}$ above it. With some operations $B^r$ we compute $\codeg(\D) > i$, contradiction.

	\item Assume that the result holds up to $k-1$ for some $n\leq k-1 \leq N-1$. 
	All the vertices with $r(v) \leq k-1$ have been chosen and are below $v_{\alpha_{k-1}+i}$ ; and $a_k-i$ vertices with $r(v)=k$ are between $v_{\alpha_{k-1}-i+1}$ and $v_{\alpha_{k-1}+i}$. Thus it remains $i$ vertices with $r(v)=k$ to choose. 
	If $r(v_j)=k$ with $j > \alpha_k+i$, then between $v_{\alpha_k-i+1}$ and $v_{\alpha_k+i}$ there are at most $i-1$ vertices having $r(v)=k$, so at least $i+1$ vertices having $r(v) > k$. These $i+1$ vertices are all below $v_j$ thus we can perfom at least $i+1$ operations $B^r$ each of them making the codegree drop by at least $1$. Hence $\codeg(\D) > i$, contradicton.
	
	Thus between $v_{\alpha_k-i+1}$ and $v_{\alpha_k+i}$ there are $i$ vertices with $r(v) = k$. Assume that $r(v_j)>k+1$ for some $\alpha_k-i+1 \leq j \leq \alpha_k+i$. Then $v_j$ has at most $j-\alpha_k+i$ vertices with $r(v) \in \{k, k+1\}$ below it, so at least 
	\[ (i+a_{k+1}) - (j-\alpha_k+i) = \alpha_k+a_{k+1}-j \geq a_{k+1}-i > i\]
	vertices with $r(v) \in \{k, k+1\}$ above it. With some operations $B^r$ we compute $\codeg(\D) > i$, contradiction.
	\end{itemize}
	
	\item The proof is similar as (1). If $r(v_j) < N$ for some $\alpha_N-i+1 \leq j \leq \alpha_N$, then the vertex $v_j$ has $j$ vertices with $r(v)=N$ below it, thus we can perform $j$ operations $B^r$, each of them making the codegree drop by at most $1$ by lemma \ref{lemme_opeAB}. Since $j \geq \alpha_N-i+1$ and $\alpha_N \geq a_N>2i$ we get $\codeg(\D) \geq i+1$, a contradiction.
\end{enumerate}
\end{demo}

\begin{rk}\label{rk_r}
The point $(1)$ is also true if the order is total only on the $b$ lowest vertices, with some $b\geq a_n-i$ ; the proof is the same.
\end{rk}

\begin{rk}
 Consider $\Delta'$ the symmetric of $\Delta$ with respect to a vertical axis. Then any diagram $\D'$ with Newton polygon $\Delta'$ corresponds to a unique diagram $\D$ with Newton polygon $\Delta$. Their functions $r',\ell'$ and $r,\ell$ are linked by $r' = -\ell$ and $\ell' = -r$. Thus, applying lemma \ref{lemme_r} to $\Delta'$ gives a similar result for the function $\ell$.
\end{rk}

\begin{lemme} \label{lemme_rgamma}
Let $i \in \N$ and $\Delta$ be a non-singular polygon with $\Delta > 2i$. The number of possible couples $(r,\ell)$ to construct a floor diagram with Newton polygon $\Delta$,
 codegree at most $i$ and having a total order on its vertices is 
\[ \dsum_{k=0}^i \ \ \dsum_{k_1 + \dots +k_{\chi^*(\Delta)-2} = k} p(k_1)\dots p(k_{\chi^*(\Delta)-2}) \]
where $\chi^*(\Delta)$ is the number of non-horizontal edges of $\Delta$.
\end{lemme}

\begin{demo}
We use the same notations as in lemma \ref{lemme_r}.

	By lemma \ref{lemme_r}, the function $r$ is entirely described by the data of vectors $\tilde\gamma^k \in(\N^*)^i$ for any $n\leq k\leq N-1$
	such that the vertices between $v_{\alpha_k-i+1}$ and $v_{\alpha_k+i}$ having $r(v)=k$ are the $v_{\alpha_k-i+\tilde\gamma^k_j}$ for $1 \leq j \leq i$. 
	
	Given $n \leq k \leq N-1$ and $1 \leq j \leq i$, the vertex $v_{\alpha_k-i+\tilde\gamma^k_j}$ has $\tilde\gamma^k_j-j$ vertices with $r(v)=k+1$ below it, thus we can perform 
	\[ c_r := \dsum_{k=n}^{N-1} \dsum_{j=1}^i (\tilde\gamma^k_j -j) \]
	operations $B^r$, each of them making the codegree drop by $1$. Hence the function $r$ contributes $c_r$ to the codegree of $\D$.
	
	Similarly, by the previous remark the function $\ell$ is entirely determined by vectors $\tilde\delta^k \in(\N^*)^i$ for any $n'\leq k\leq N'-1$, where $(n',-1), (n'+1,-1), \dots, (N',-1)$ are the primitive direction vectors of the edges of the left side of $\Delta$, from top to bottom. Similarly, we can perform 
	\[ c_\ell := \dsum_{k=n'}^{N'-1} \dsum_{j=1}^i (\tilde\delta^k_j -j) \]
	operations $B^\ell$, each of them making the codegree drop by $1$. Hence the function $\ell$ contributes $c_\ell$ to the codegree of $\D$.
	
	Putting together these two contributions we should have
	\[  c_r + c_\ell  \leq \codeg(\D) \leq i .\] 
	
	Given a vector $\tilde\beta \in (\N^*)^i$, we consider the vector $\beta \in \N^i$ defined by $\beta_j = \tilde\beta_{i-j+1} - \tilde\beta_{i-j} -1$, where $\tilde\beta_0=0$ by convention. 
	One has	
	\[ \codeg(\beta) = \dsum_{j=1}^i (\tilde\beta_j -j) . \]
	Applying this to the vectors $(\tilde\gamma^k)_k$ and $(\tilde\delta^k)_k$, we see that the data of functions $r$ and $\ell$ satisfying $c_r + c_\ell \leq i$ is equivalent to the data of vectors $(\gamma^k)_k$ and $(\delta^k)_k$ satisfying
	\[ \dsum_{k=n}^{N-1} \codeg(\gamma^k) + \dsum_{k=n'}^{N'-1} \codeg(\delta^k) \leq i ,  \]
\ie to the data of $N+N'-n-n' = \chi^*(\Delta)-2$ vectors whose sum of codegrees is at most $i$. If that sum equals $k$, this corresponds to the data of a decomposition $k_1 + \dots + k_{\chi^*(\Delta)-2} = k$, and then for any $1 \leq j \leq  \chi^*(\Delta)-2$ to the data of a vector of codegree $k_j$. By lemma \ref{lemme_cardBi} there are $p(k_j)$ possibilities for such a vector. Hence, we conclude that the number of possible couples $(r,\ell)$ is
\[  \dsum_{k=0}^i \ \ \dsum_{k_1 + \dots +k_{\chi^*(\Delta)-2} = k} p(k_1)\dots p(k_{\chi^*(\Delta)-2}) . \]
\end{demo}

\begin{theo} \label{theo_nonsinghori}
Let $i \in \N$ and $\F$ be a $h$-transverse and non-singular fan having rays generated by $(0,1)$ and $(0,-1)$. Let $\Delta \in D(\F)$ and $s\in \{0,\dots, \smax(\Delta) \}$. If $(\Delta,s)$ satisfies 
\[ \left\{ \begin{array}{rcl}
 \Delta &>& 2(i+2) \\
 e^{-\infty}(\Delta) &>& i+2s \\
 e^{\pm \infty}(\Delta) &>& i+ (a(\Delta)-\midaD+1)d_\F
\end{array} \right. \]
then 
\[ \< \GDelta(s)\>_i = P_i(y(\Delta),\chi(\Delta),s). \]
\end{theo}

\begin{demo}
Let $i$, $\F$, $\Delta$ and $s$ be as in the hypothesis.
We look for a formula for $\< \GDelta(s) \>_i$. We  will use the pairing $\{\{1,2\},\dots,\{2s-1,2s  \}\}$ of order $s$. We use the shortcuts $(a,e^{+\infty},e^{-\infty},\chi) = (a(\Delta),e^{+\infty}(\Delta),e^{-\infty}(\Delta),\chi(\Delta))$.

	By lemma \ref{lemmetotalorder} the order is total on the vertices of any diagram $\D$ having Newton polygon $\Delta$ and codegree at most $i$ ; we denote them by $v_1 \prec \dots \prec v_a$. Any diagram $\D$ has the shape of figure \ref{overallshape} : it can be entirely described by the data of the vectors $u,\utilde \in \N^i$ and of the functions $r,\ell$, hence we will use the notation $\D = (u,\utilde,r,\ell)$.
	
\begin{figure}[h]
	\centering
	\begin{tikzpicture}[scale=8/9]
	
	\floor (1) at (0,1.5) {\scriptsize $\ell(v_1) \ \ \ \ \ r(v_1)$} ;
	\floor (2) at (1.5,3) {\scriptsize $\ell(v_2) \ \ \ \ \ r(v_2)$} ;
	\floor (3) at (3.5,5) {\scriptsize $\ell(v_{i+1}) \ \ \ \ \ r(v_{i+1})$} ;
	\floor (4) at (3.5,7) {\scriptsize $\ell(v_{a-i}) \ \ \ \ \ r(v_{a-i})$} ;
	\floor (5) at (1.5,9) {\scriptsize $\ell(v_{a-1}) \ \ \ \ \ r(v_{a-1})$} ;
	\floor (6) at (0,10.5) {\scriptsize $\ell(v_a) \ \ \ \ \ r(v_a)$} ;
	
	\draw (1) to (2) ;
	\draw[dashed] (2) to (3) ;
	\draw[dashed] (3) to (4) ;
	\draw[dashed] (4) to (5) ;
	\draw (5) to (6) ;

		\draw (1) to[out=-120, in=90] (-0.4,0) ;
		\draw (1) to[out=-60, in=90]  (0.4,0);
		\node at (0,0.6) {$\cdots$} ;
		\node at (0,-0.6) {};

		\draw (2) to[out=-70, in=90] (2.2,0) ;
		\draw (2) to[out=-40, in=90]  (3,0);
		
		\node at (2.6,0.6) {$\cdots$} ;
		\node at (2.6,-0.4) {\large $\underbrace{\ \ \ \ \ }_{u_1}$} ;
		
		\node at (3.6,0.6) {$\cdots$} ;
		
		\draw (3) to[out=-75, in=90] (4.2,0) ;
		\draw (3) to[out=-50, in=90]  (5,0);
		\node at (4.6,0.6) {$\cdots$} ;
		\node at (4.6,-0.4) {\large $\underbrace{\ \ \ \ \ }_{u_i}$} ;

		\draw (6) to[out=120, in=-90] (-0.4,12) ;
		\draw (6) to[out=60, in=-90]  (0.4,12);
		\node at (0,11.4) {$\cdots$} ;
		\node at (0,12.6) {};
		
		\draw (5) to[out=70, in=-90] (2.2,12) ;
		\draw (5) to[out=40, in=-90]  (3,12);
		\node at (2.6,11.4) {$\cdots$} ;
		\node at (2.6,12.4) {\large $\overbrace{\ \ \ \ \ }^{\utilde_1}$} ;
		
		\node at (3.6,11.4) {$\cdots$} ;
		
		\draw (4) to[out=75, in=-90] (4.2,12) ;
		\draw (4) to[out=50, in=-90]  (5,12);
		\node at (4.6,11.4) {$\cdots$} ;
		\node at (4.6,12.4) {\large $\overbrace{\ \ \ \ \ }^{\utilde_i}$} ;
	\end{tikzpicture}
	\caption{Overall shape for $\D$.}
	\label{overallshape}
\end{figure}
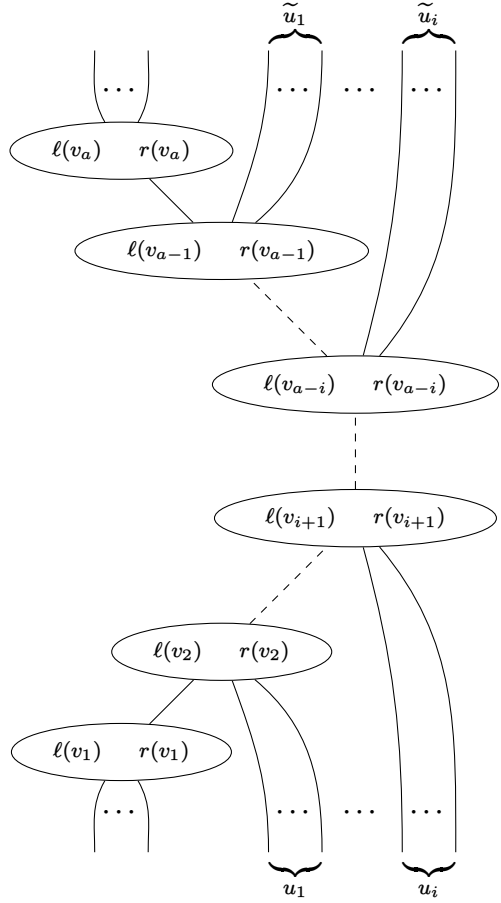

	Since $e^{-\infty} > i+2s$ then for any compatible marking $m$ and any $j \leq 2s$ we have $m^{-1}(j) \in E^{-\infty}(\D)$. In particular, the multiplicity of $(\D,m)$ does not depend on $m$ and is just
	\[ \mu(\D) = \dprod_{e\in E_0} [w(e)]^2  \]
We then have	
	\[ \< \GDelta(s) \>_i = \dsum_{\D} \nu(\D) \< \mu(\D)\>_{i-\codeg(\D)}  \]
	where the sum runs over the floor diagrams of Newton polygon $\Delta$ 
	 and codegree at most $i$, and where $\nu(\D)$ is the number of compatible markings of $\D$.
	For $\D = (u,\utilde,r,\ell)$ this number only depends on $u$ and $\utilde$ and is
	\[ \nu(\D) =  \dsum_{S\in d(s)} \binom{s}{S} \nu_{\geq1}(e^{-\infty},2,u,S) \nu_{\geq1}(e^{+\infty},2,\utilde,0). \]	
	By lemma \ref{lemmepoidsgrand} and lemma \ref{lemme_mult} its multiplicity gives
	\[ \< \mu(\D) \>_{i-\codeg(\D)} = \Phi_{i-\codeg(\D)}(a-1) \]
	which is also independent of $r$ and $\ell$. Thus, to compute $\< \GDelta(s)\>_i$ we need to determine how many couples $(r,\ell)$ are possible and then sum over $(u,\utilde)$. 
	Here, $\chi^*(\Delta) = \chi-2$ so by lemma \ref{lemme_rgamma} the number of possible couples $(r,\ell)$ is
\[  \dsum_{k=0}^i \ \ \dsum_{k_1 + \dots +k_{\chi-4} = k} p(k_1)\dots p(k_{\chi-4}) . \]

	Once $(r,\ell)$ is chosen and contributes $k$ to the codegree, it remains to sum over $(u,\utilde)$ such that $\codeg(u+\utilde) \leq i-k$, \ie $u+\utilde \in C_{i-k}$. We can now compute 
\begin{align*}
\< \GDelta(s) \>_i &= \dsum_{\D} \nu(\D) \< \mu(\D) \>_{i-\codeg(\D)} \\
  &= \dsum_{k=0}^i \ \  \dsum_{k_1 + \dots +k_{\chi-4} = k} p(k_1)\dots p(k_{\chi-4}) \dsum_{u+\utilde \in C_{i-k}} \nu(\D) \Phi_{i-k-\codeg(u+\utilde)}(a-1) \\
  &= \dsum_{k=0}^i \ \ \dsum_{k_1 + \dots +k_{\chi-4} = k} p(k_1)\dots p(k_{\chi-4}) \dsum_{j=1}^{i-k} \ \dsum_{u+\utilde \in B_j} \nu(\D) \Phi_{i-k-j}(a-1) 
\end{align*}
	which shows that $\< \GDelta(s) \>_i$ coincides with the degree $i$ coefficient of the product of : 
\begin{itemize}
	\item the generating series of $\left(\dsum_{k_1 + \dots +k_{\chi-4} = k} p(k_1)\dots p(k_{\chi-4})\right)_k$, which is $A_2^{\chi-4}$ by definition,
	
	\item the generating series of $\left(\dsum_{u+\utilde \in B_k} \nu(\D)\right)_k$, which is $A_0^s A_1^{e^{-\infty}+e^{+\infty}-2s}A_2^4$ by lemma \ref{lemme_binom},
	
	\item the generating series of $(\Phi_k(a-1))_k$, which is $A_1^{2a-2}$ by lemma \ref{coro_genPhi}.
\end{itemize}
	Since $y= e^{-\infty}+e^{+\infty} +2a$, this product is 
	\[  A_0^s A_1^{y-2-2s} A_2^\chi \]
and its degree $i$ coefficient is a polynomial $P_i$ of degree $i$ in the variables $y$, $\chi$ and $s$.
\end{demo}

\subsection{The blow-up trick and the case of $\CP^2$}
\label{subsec_blowupCP2}

Let $\F$ be a fan. Let $u$ and $v$ be two primitive generators of two consecutive rays of $\F$. Consider $\tilde\F$ be the fan constructed from $\F$ by adding a ray generated by $u+v$. We say that $\tilde\F$ is a \emph{blow-up} of $\F$. This terminology comes from toric geometry : the surface $X_{\tilde\F}$ is the blow-up of $X_\F$ in one point. At the level of the dual polygons, any $\tilde\Delta \in D(\tilde\F)$ is obtained by cutting off a corner of a $\Delta \in D(\F)$, and is said to be a blow-up of $\Delta$.

Let $i\in\N$. For $\Delta$ a polygon with dual fan $\F$ and $s$ an integer we consider the following conditions :
\[ \stari \left\{ \begin{array}{rcl}
 \Delta &>& 2(i+2) \\
 e^{-\infty}(\Delta) &>& i+2s \\
 e^{-\infty}(\Delta) &>& i+ (a(\Delta)-\midaD+1)d_\F
\end{array} \right. . \] 

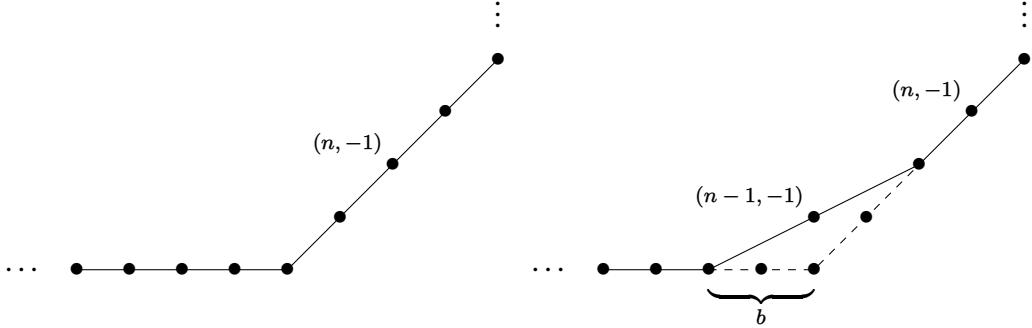
\begin{figure}[h] 
		\centering
	\begin{tikzpicture}[scale=0.7]
		
		\draw (-4,0) to (0,0) to node[above left] {\scriptsize $(n,-1)$} (4,4) ;
		\foreach \p in {(-4,0),(-3,0),(-2,0),(-1,0), (0,0), (1,1),(2,2),(3,3),(4,4)}
			{\draw node at \p {$\bullet$} ;} 
		\draw node at (-5,0) {$\dots$} ;
		\draw node at (4,5) {$\vdots$} ;

		\draw (6,0) to (8,0) to node[above left] {\scriptsize $(n-1,-1)$} (12,2) to node[above left] {\scriptsize $(n,-1)$} (14,4);
		\draw[dashed] (8,0) to (10,0) to (12,2) ;
		\foreach \p in {(6,0),(7,0),(8,0),(9,0), (10,0),(10,1),(11,1),(12,2),(13,3),(14,4)}
			{\draw node at \p {$\bullet$} ;} 
		\draw node at (5,0) {$\dots$} ;
		\draw node at (14,5) {$\vdots$} ;
		
		\draw node at (9,-0.7) {$\underbrace{\ \ \ \ \ \ \ \ \ \ \ }_{b}$} ;

	\end{tikzpicture}
	\caption{The polygons $\Delta$ (left) and $\tilde\Delta$ (right).} 
	\label{fig_blowup} 
\end{figure}

\begin{prop} \label{prop_blowup}
 Let $\Delta$ and $\tilde\Delta$ be the polygons whose bottom right corners are depicted in figure \ref{fig_blowup}. 
 There exists a series $B$ such that if $(\Delta,s)$ satisfies $\stari$ then $\< \GDelta(s) \>_i$ is given by the degree $i$ coefficient of 
 \[ A_0^s A_1^{e^{-\infty}(\Delta)-2s} A_2^2 \times B ,\]
 and if $(\tilde\Delta,s)$ satisfies $\stari$ then $\< G_{\tilde\Delta}(0;s) \>_i$ is given by the degree $i$ coefficient of 
 \[ A_0^s A_1^{e^{-\infty}(\tilde\Delta)-2s} A_2^3 \times B .\]
\end{prop}
  
\begin{demo}
	 Let $a=a(\Delta) = a(\tilde\Delta)$ and $e^{-\infty}=e^{-\infty}(\Delta)$. We denote by $a_n$ the integer length of the edge of the right side of $\Delta$ having direction vector $(n,-1)$. We denote by $b$ the integer length of a side of the triangle we cut off to obtain $\tilde\Delta$, see figure \ref{fig_blowup}. 
	We will use the pairing $\{\{1,2\},\dots,\{2s-1,2s  \}\}$ of order $s$.  
	
	We will first make a calculation for $\Delta$. Then we will explain how to construct a correspondence between the floor diagrams with Newton polygon $\Delta$ and the ones with Newton polygon $\tilde\Delta$, allowing us to make a calculation for $\tilde\Delta$.

\tocless{\subsection*{Calculation for $\Delta$.}}
Because $\Delta$ is large enough, any diagram $\D$ that contributes to $\< \GDelta(s)\>_i$ has a unique minimal vertex by lemma \ref{lemmetotalorder}. By lemma \ref{lemme_ordertop}, $\D$ has a total order on its $a-i-1$ lowest vertices. It may have several maximal vertices that we hide in a very top part $\T$, see figure \ref{fig_cappingtreebis}.
 By remark \ref{rk_r}, the function $r$ is constant on the $a_n-i$ lowest vertices. 
Moreover, we have the following inequalities :
\[ i+1 < b+i < \Midab < a-i-1  .\]
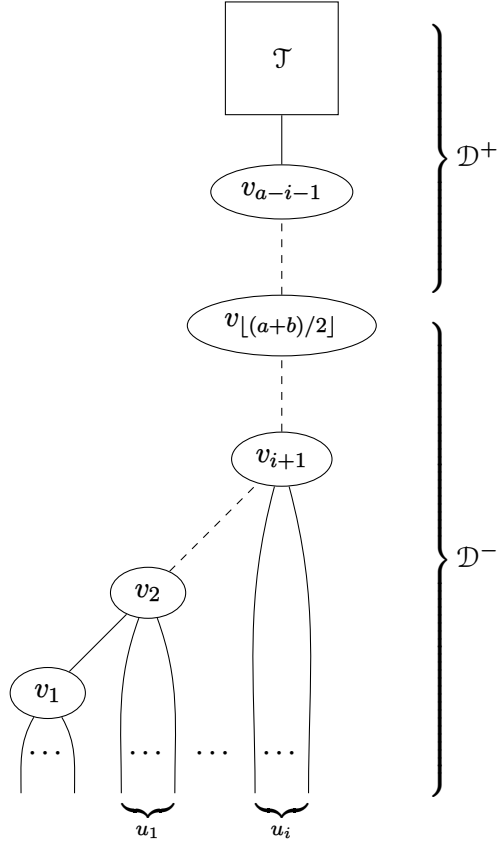
\begin{figure}[h!]
	\centering
	\begin{tikzpicture}[scale=8/9]
	
	\floor (1) at (0,1.5) {$v_1$} ;
	\floor (2) at (1.5,3) {$v_2$} ;
	\floor (3) at (3.5,5) {$v_{i+1}$} ;
	\floor (4) at (3.5,7) {$v_{\midab}$} ; 
	\floor (5) at (3.5,9) {$v_{a-i-1}$} ;
	\node[draw,rectangle, minimum width=1.5cm, minimum height=1.5cm] (6) at (3.5,11) {$\T$} ;
	
	\draw (1) to (2) ;
	\draw[dashed] (2) to (3) ;
	\draw[dashed] (3) to (4) ;
	\draw[dashed] (4) to (5) ;
	\draw (5) to (6) ;

	\draw (1) to[out=-120, in=90] (-0.4,0) ;
	\draw (1) to[out=-60, in=90]  (0.4,0);
	\node at (0,0.6) {$\cdots$} ;
	\node at (0,-0.6) {} ;

	\draw (2) to[out=-110, in=90] (1.1,0) ;
	\draw (2) to[out=-70, in=90]  (1.9,0);
	\node at (1.5,0.6) {$\cdots$} ;
	\node at (1.5,-0.4) {\large $\underbrace{\ \ \ \ \ }_{u_1} $} ;
		
	\node at (2.5,0.6) {$\cdots$} ;
		
	\draw (3) to[out=-105, in=90] (3.1,0) ;
	\draw (3) to[out=-75, in=90]  (3.9,0);
	\node at (3.5,0.6) {$\cdots$} ;
	\node at (3.5,-0.4) {\large $\underbrace{\ \ \ \ \ }_{u_i}$} ;
	\node at (6,3.5) {$ \left. \begin{array}{c}
 \\ \\ \\ \\	 \\ \\ \\ \\	 \\ \\ \\ \\	 \\ \\  
\end{array}	 \right\} \D^-$ } ;
	\node at (6,9.5) {$ \left. \begin{array}{c}
 \\ \\ \\ \\	 \\ \\ \\ \\	 
\end{array}	 \right\} \D^+$ } ;
	
	\end{tikzpicture}
	\caption{Decomposition of a diagram.}
	\label{fig_cappingtreebis}
\end{figure}
 We can cut $\D$ into two parts, see figure \ref{fig_cappingtreebis} :
 \begin{itemize}
  \item a bottom part : we denote by $\D^-$ the diagram consisting of the sources, the vertices from $v_1$ to $v_\midab$ and the bounded edges between them ; 
  \item a top part : we denote by $\D^+$ the remaining of $\D$ ; it has $a-\midab$ vertices.
 \end{itemize}
This leads to consider the following sets. Let $k\in \{0,\dots,i\}$. We define
\begin{itemize}
	\item $B_k(\Delta)$ the set of all possible bottom parts $\D^-$ having $\midab$ vertices and of codegree $k$. Encoding how the infinite edges oriented inward are attached to the vertices we establish a bijection between $B_k(\Delta)$ and $B_k$ : each $\D^-$ can be represented by a $u\in B_k$.
	
	\item $T_{i-k}(\Delta)$ the set of all possible top parts $\D^+$ having $a-\midab$ vertices and of codegree at most $i-k$.
\end{itemize}
As explained above there is a bijection
\[ C_i(\Delta) \simeq \bigsqcup_{k=0}^i B_k \times T_{i-k}(\Delta).  \]
The number of markings of a diagram $\D$ can be calculated separately on its top and bottom parts. If $\D$ is represented by $(u,\D^+) \in B_k \times T_{i-k}(\Delta)$, we denote by $\nu(\D^+)$ the number of markings of the top part $\D^+$, and the number of markings of the bottom part is
\[ \dsum_{S \in d(s)} \binom{s}{S} \nu_{\geq1}(e^{-\infty}-2s,2,u,S).\]
Moreover, because $e^{-\infty} > i+2s$ then for any compatible marking $m$ and any $j\leq 2s$ one has $m^{-1}(j) \in E^{-\infty}(\D)$. Hence the multiplicity of $(\D,m)$ does not depend on $m$ and is
\[ \mu(u,\D^+) =  \mu(\D) = \dprod_{e\in E_0} [w(e)]^2 = \dprod_{e\in E_0 \cap E^(\D^-)} [w(e)]^2 \times \dprod_{e\in E_0 \cap E^(\D^+)} [w(e)]^2   = \mu(\D^-) \times \mu(\D^+) .\]
By lemmas \ref{lemmepoidsgrand} and \ref{lemme_mult} it gives
\begin{align*}
\< \mu(\D) \>_{i-k-\codeg(\D^+)}
  &= \dsum_{i_1+i_2=i-k-\codeg(\D^+)} \< \mu(\D^+)\>_{i_1} \< \mu(\D^-)\>_{i_2}   \\
  &= \dsum_{i_1+i_2=i-k-\codeg(\D^+)} \< \mu(\D^+)\>_{i_1}  \Phi_{i_2}(\midab-1) 
\end{align*}
and this is independent of $u$.
 Hence if we set 
\begin{align*}
\alpha_k(e^{-\infty},s) &= \dsum_{u \in B_k} \dsum_{S \in d(s)} \binom{s}{S} \nu_{\geq1}(e^{-\infty}-2s,u,2,S), \\
\beta_k(a,b) &= \dsum_{\D^+ \in T_k(\Delta)} \nu(\D^+) \dsum_{i_1+i_2=k-\codeg(\D^+)} \< \mu(\D^+)\>_{i_1}  \Phi_{i_2}(\midab-1),
\end{align*}
then we get 
\begin{align*}
\< G_{\Delta}(0;s) \>_i
  &= \dsum_{\codeg(\D)\leq i} \nu(\D) \< \mu(\D) \>_{i-\codeg(\D)} \\
  &= \dsum_{k=0}^i \dsum_{u\in B_k} \dsum_{S \in d(s)} \binom{s}{S} \nu_{\geq1}(e^{-\infty}-2s,2,u,S) \dsum_{\D^+ \in T_{i-k}(\Delta)}  \nu(\D^+) \< \mu(u,\D^+) \>_{i-k-\codeg(\D^+)} \\
  &= \dsum_{k=0}^i \alpha_k(e^{-\infty},s)\beta_{i-k}(a,b)
\end{align*}
which shows that $\< G_{\Delta}(0;s) \>_i$ is the degree $i$ coefficient of the product of
\begin{itemize}
	\item the generating series of $(\alpha_k(e^{-\infty},s))_k$ which is $A_0^s A_1^{e^{-\infty}-2s} A_2^2$ by lemma \ref{lemme_binom},
	
	\item the generating series of $(\beta_k(a,b))_k$ that we will denote by $B$.
\end{itemize}

\tocless{\subsection*{Calculation for $\tilde\Delta$.}}
Now let's have a look at what happens for $\tilde\Delta$. Let $k \in \{0,\dots,i\}$ and $j\in\{0,\dots,i-k\}$ be integers, let $\gamma\in B_k$ be a vector of codegree $k$, and let $\D = (u,\D^+)\in B_j\times T_{i-k-j}$ be a diagram with Newton polygon $\Delta$, genus $0$ and codegree at most $i-k$. Because $b>2(i+2)$ we can see $\gamma$ as a vector in $\N^b$. 
Consider $\tilde\gamma \in\N^b$ the vector given by
\[ \tilde\gamma_m = \dsum_{\ell=0}^{m-1} (\gamma_{b-\ell}+1),\ \text{ \ie }\ \gamma_m = \tilde\gamma_{b-m+1} - \tilde\gamma_{b-m} -1 \]
with the convention $\tilde\gamma_0=0$. Note that we have
$\tilde\gamma_1 \geq 1$, $\tilde\gamma_m > \tilde\gamma_{m-1}$ and
\[ \codeg(\gamma) = \dsum_{m=1}^b (\tilde\gamma_m-m). \]
In particular, $\tilde\gamma_m \leq b+i$ otherwise we would have $\codeg(\gamma) > k$. This implies $\tilde\gamma_m \leq a_n-i$ and so $r(v_{\tilde\gamma_m}) = n$ by lemma \ref{lemme_r}.

Consider the diagram $\D_\gamma$ obtained from $\D$ with the following process : set $r(v_{\tilde\gamma_m})=n-1$ for any $1 \leq m \leq b$, adjust the weight of the bounded edges and remove enough sources to the minimal vertex of $\D$ to satisfy the divergence condition. Then $\D_\gamma \in C_i(\tilde\Delta)$ ; its codegree is $\codeg(\D_\gamma) = \codeg(\D) + \codeg(\gamma) \leq i$. Note that we did not change anything in the top part $\D^+$ of $\D$.

	Conversely, let $\tilde\D \in C_i(\tilde\Delta)$. By lemmas \ref{lemmetotalorder} and \ref{lemme_ordertop}, $\tilde\D$ admits a total order on its $a-i-1$ lowest vertices.
	The diagram $\tilde\D$ has $b$ vertices with $r(v)=n-1$. Suppose there is such a vertex with more than $b+i$ vertices below it. Then at least $i+1$ of the vertices below it have $r(v)=n$. Performing $i+1$ operations $B^r$ we see that $\codeg(\tilde\D) \geq i+1$, a contradiction. Thus all the vertices with $r(v) = n-1$ are between $v_1$ and $v_{b+i}$.
	We denote by $\tilde\gamma \in \N^b$ the vector whose coordinates are the indices of the vertices with $r(v) = n-1$, and let $\gamma \in\N^b$ be defined by $\gamma_m = \tilde\gamma_{b-m+1} - \tilde\gamma_{b-m} -1$.
	Then 
	\[ k := \codeg(\gamma) = \dsum_{m=1}^b (\tilde\gamma_m-m) \leq \codeg(\tilde\D) \leq i \]
	and $\tilde\D = \D_\gamma$, where $\D$ is the diagram of $C_{i-k}(\Delta)$ obtained by setting $r(v_{\tilde\gamma_m})=n$ for any $1 \leq m \leq b$, then adjusting the weights and adding enough sources to $v_1$.

	In other words, there is a bijection
\[ C_i(\tilde\Delta) \simeq \bigsqcup_{k=0}^i B_k \times C_{i-k}(\Delta) \]
and a floor diagram with Newton polygon $\tilde\Delta$, genus $0$ and codegree at most $i$ can be represented as
\[ \tilde\D = \D_\gamma = (\gamma,\D) \in B_k \times C_{i-k}(\Delta) \]
for some $1\leq k \leq i$. 
Its codegree is 
\[ \codeg(\tilde\D) = \codeg(\gamma) + \codeg(\D) = k+\codeg(\D) \]
The diagram $\D$ is itself represented by
\[ \D = (u,\D^+) \in B_j \times T_{i-k-j}(\Delta) \]
for some $1\leq j \leq i-k$, so 
\[ \tilde\D = (\gamma,u,\D^+) \in B_k \times B_j \times T_{i-k-j}(\Delta). \]
 The number of markings can be computed separately on a bottom part $\tilde\D^- = (\gamma,u)$ and on a top part $\tilde\D^+ = \D^+$.
	Moreover, by lemma \ref{lemmepoidsgrand} all the edges of $\tilde\D^-$ have a weight greater than $i-\codeg(\tilde\D)$  so that in the expression
\begin{multline*}
   \< G_{\tilde\Delta}(0;s) \>_i 
  = \dsum_{k=0}^i \dsum_{\gamma \in B_k} \dsum_{j=0}^{i-k} \dsum_{u\in B_j} \dsum_{S\in d(s)} \binom{s}{S} \nu_{\geq1}(e^{-\infty}-b-2s,2,u,S) \\
  \dsum_{\D^+ \in T_{i-k-j}(\Delta)} \nu(\D^+) \< \mu(\gamma,u,\D^+) \>_{i-k-j-\codeg(\D^+)} %
\end{multline*}
the inner sum does not depend neither on $u$ nore on $\gamma$ by lemma \ref{lemme_mult}. Hence, the sum over $\gamma \in B_k$ contributes $p(k)$ by lemma \ref{lemme_cardBi}, and we get
\begin{align*}
 \< G_{\tilde\Delta}(0;s) \>_i
  &= \dsum_{k=0}^i p(k) \dsum_{j=0}^{i-k} \alpha_j(e^{-\infty}-b,s) \beta_{i-k-j}(a,b)
\end{align*}
which shows that $\< G_{\tilde\Delta}(0;s) \>_i$ is the degree $i$ coefficient of the product of 
\begin{itemize}
	\item the generating series of $(p(k))_k$ which is $A_2$ by definition, 
	
	\item the generating series of $(\alpha_k(e^{-\infty}-b,s))_k$ which is $A_0^s A_1^{e^{-\infty}-b-2s} A_2^2$ by lemma \ref{lemme_binom},
	
	\item the generating series of $(\beta_k(a,b))_k$ which is $B$.
\end{itemize}
Since $e^{-\infty}(\tilde\Delta) = e^{-\infty}-b$ we can conclude.
\end{demo}

\begin{coro} \label{coro_blowup}
With the same notations as in proposition \ref{prop_blowup}, let $i_m$ be an integer such that both $(\Delta,s)$ and $(\tilde\Delta,s)$ satisfy $(\star)_{i_m}$. The following are equivalent.
\begin{enumerate}
\item For any $i\in\{0,\dots,i_m\}$, $ \< \GDelta(s) \>_i = P_i(y(\Delta),\chi(\Delta),s)$,
 
\item For any $i\in\{0,\dots,i_m\}$, $ \< G_{\tilde\Delta}(0;s) \>_i = P_i(y(\tilde\Delta),\chi(\tilde\Delta),s)$.
\end{enumerate}
\end{coro}

\begin{demo}
For any $0 \leq i \leq i_m$, both $(\Delta,s)$ and $(\tilde\Delta,s)$ satisfy $\stari$.
 Hence by definition of $(P_i)_i$ and proposition \ref{prop_blowup} one has the following equivalences :
 \begin{align*}
 (1) &\Leftrightarrow A_0^s A_1^{e^{-\infty}(\Delta)-2s} A_2^2 \times B = A_0^s A_1^{y(\Delta)-2-2s} A_2^{\chi(\Delta)} \mod x^{i_m} \\
  &\Leftrightarrow A_0^s A_1^{e^{-\infty}(\tilde\Delta)+b-2s} A_2^3 \times B = A_0^s A_1^{y(\tilde\Delta)+b-2-2s} A_2^{\chi(\Delta)+1} \mod x^{i_m} \\
  &\Leftrightarrow  A_0^s A_1^{e^{-\infty}(\tilde\Delta)-2s} A_2^3 \times B = A_0^s A_1^{y(\tilde\Delta)-2-2s} A_2^{\chi(\tilde\Delta)} \mod x^{i_m}  \\
  &\Leftrightarrow (2) . 
 \end{align*}
\end{demo}

This can be used to compute $\< \GDelta(0,s) \>_i$ when the underlying toric surface is $\CP^2$.

\begin{coro} \label{coro_CP2}
Let $\F$ be the fan whose rays are generated by $(-1,0)$, $(0,-1)$ and $(1,1)$. For any $i\in\N$, any $\Delta \in D(\F)$ and $s\in\{0,\dots,\smax(\Delta) \}$, if 
\[  \left\{ \begin{array}{rcl}
	\Delta &>& 5(i+1)+6 \\
	\Delta &>& i+2s
\end{array} \right. \]
 then one has 
\[ \< \GDelta(s) \>_i = P_i(y(\Delta),3,s).\] 
\end{coro}

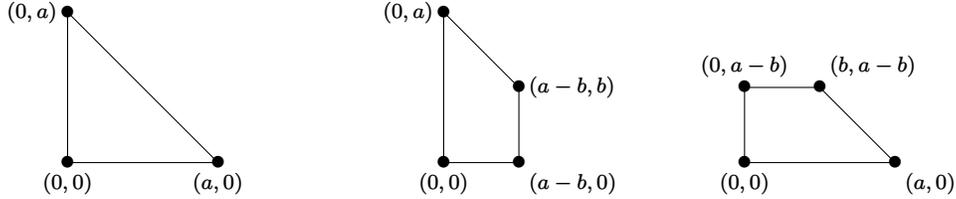
\begin{figure}[h] 
	\centering
	\begin{tikzpicture}[scale=1] 
		\draw (0,0) node {$\bullet$} ;
		\draw (0,0) node[below] {\scriptsize $(0,0)$} ;
		\draw (2,0) node {$\bullet$} ;
		\draw (2,0) node[below] {\scriptsize $(a,0)$} ;
		\draw (0,2) node {$\bullet$} ; 
		\draw (0,2) node[left] {\scriptsize $(0,a)$} ; 
		\draw (0,0) -- (2,0) ;
		\draw (0,0) -- (0,2) ;
		\draw (2,0) -- (0,2) ;
		
		\draw (5,0) node {$\bullet$} ;
		\draw (5,0) node[below] {\scriptsize $(0,0)$} ;
		\draw (6,0) node {$\bullet$} ;
		\draw (6,0) node[below right] {\scriptsize $(a-b,0)$} ;
		\draw (5,2) node {$\bullet$} ; 
		\draw (5,2) node[left] {\scriptsize $(0,a)$} ; 
		\draw (6,1) node {$\bullet$} ; 
		\draw (6,1) node[right] {\scriptsize $(a-b,b)$} ; 
		\draw (5,0) -- (6,0) -- (6,1) -- (5,2) -- (5,0) ;
		
		\draw (9,0) node {$\bullet$} ;
		\draw (9,0) node[below] {\scriptsize $(0,0)$} ;
		\draw (11,0) node {$\bullet$} ;
		\draw (11,0) node[below right] {\scriptsize $(a,0)$} ;
		\draw (10,1) node {$\bullet$} ; 
		\draw (10,1) node[above right] {\scriptsize $(b,a-b)$} ; 
		\draw (9,1) node {$\bullet$} ; 
		\draw (9,1) node[above] {\scriptsize $(0,a-b)$} ; 
		\draw (9,0) -- (11,0) -- (10,1) -- (9,1) -- cycle ;
	\end{tikzpicture}
	\caption{The polygons $\Delta_a$ (left) and $\Delta_{a,b}$ (middle) and $\Delta'_{a,b}$ (right).}
	\label{fig_blowupCP2}
\end{figure}

\begin{demo}
A polygon $\Delta \in D(\F)$ is a $\Delta_a$ for some $a\in\N^*$, see figure \ref{fig_blowupCP2}. Let $i_m$ be the maximal integer such that 
\[  \left\{ \begin{array}{rcl}
	\Delta &>& 5(i_m+1)+6 \\
	\Delta &>& i_m+2s
\end{array} \right.  \]
and $b$ be an integer such that
\[ \dfrac{a+2(i_m+1)+1}{5} < b < \dfrac{a-2(i_m+1)-1}{2}.  \]
We consider the polygon $\Delta_{a,b}$ of figure \ref{fig_blowupCP2}. 
The hypotheses imply that for any $i \in \{1,\dots,i_m\}$ one has
\[  \left\{ \begin{array}{rcl}
	a, b, a-b &>& 2(i+2) \\
	a, a-b    &>& i+2s \\
	a, a-b    &>& i + a-\Mida +1
\end{array} \right. \]
so $\Delta_a$ and $\Delta_{a,b}$ both satisfy the condition $\stari$ and we can apply corollary \ref{coro_blowup} with $\Delta = \Delta_a$ and $\tilde\Delta = \Delta_{a,b}$. We will check that the assertion $(2)$ of corollary \ref{coro_blowup} is true.

The polygon $\Delta_{a,b}$ is congruent to $\Delta'_{a,b}$, see figure \ref{fig_blowupCP2}, so that $G_{\Delta_{a,b}}(0;s) = G_{\Delta'_{a,b}}(0;s)$. Note also that $y(\Delta_{a,b}) = y(\Delta'_{a,b})$ and $\chi(\Delta_{a,b}) = \chi(\Delta'_{a,b}) = 4$.
The hypotheses imply that for any $i \in \{1,\dots,i_m\}$ one has
\[  \left\{ \begin{array}{rcl}
	a, b, a-b &>& 2(i+2) \\
	a         &>& i+2s \\
	a, b      &>& i + a-b - \Midab +1
\end{array} \right. \]
hence by theorem \ref{theo_nonsinghori} applied to $\Delta'_{a,b}$ one has  
\[ \< G_{\Delta'_{a,b}}(0;s) \>_i = P_i(y(\Delta'_{a,b}),\chi(\Delta'_{a,b}),s) = P_i(y(\Delta'_{a,b}),4,s) \]
and so $\< G_{\Delta_{a,b}}(0;s) \>_i = P_i(y(\Delta_{a,b}),4,s)$ for any $i \in \{1,\dots,i_m\}$. Hence the point $(2)$ of corollary \ref{coro_blowup} is true and we deduce that $(1)$ is also true. In particular, one has
\[ \< G_{\Delta_a}(0;s) \>_i = P_i(y(\Delta_a),\chi(\Delta),s) = P_i(y(\Delta_a),3,s) \]
for any $i\in\{0,\dots,i_m\}$.
\end{demo}

\subsection{The non-singular case}\label{subsec_maintheogen}

We can generalize theorem \ref{theo_nonsinghori} and remove the hypothesis on the existence of rays generated by $(0,1)$ and $(0,-1)$.  

\begin{lemme}\label{lemme_fan}
Let $\F$ be a non-singular and $h$-transverse fan. One of the following is true.
\begin{enumerate}
	\item The fan $\F$ is congruent to a $h$-transverse fan having two vertical rays. 
	
	\item The fan $\F$ has a single vertical ray. In that case $\F$ is obtained from a non-singular and $h$-transverse fan with 3 rays performing several blow-ups.	
\end{enumerate} 
\end{lemme}

\begin{demo}
Assume first that $\F$ has no ray generated by $(0,1)$ or $(0,-1)$. Let $n\in\N$ be the biggest integer such that $\F$ has a ray generated by $(\pm 1,\pm n)$. By symmetries we can only consider the case of $(1,n)$. Necessarily $n>0$ otherwise $\F$ is not complete. 

Let $u_1,\dots,u_k$ be the primitive generators of the rays of $\F$, taken in the anticlockwise direction and with $u_1 = (1,n)$. Write $u_2 = (a,b)$. Because $\F$ is $h$-transverse one has $a\in\{0,1,-1\}$. By hypothesis $a\neq 0$, and $a\neq 1$ otherwise $b$ would be greater than $n$. Hence $u_2 = (-1,b)$ and one has $b+n = 1$ because $\F$ is non-singular, so $u_2 = (-1,1-n)$. 

Write $u_3 = (c,d)$. By hypothesis one has $c = \pm1$ and $-d \pm (n-1) = 1$, hence $u_3=(-1,-n)$ or $u_3=(1,n-2)$. If  $u_3=(-1,-n)=-u_1$ then with the same argument one has $u_4 = -u_2$ and there is no more ray, \ie $k=4$. If  $u_3=(1,n-2)$ then the only possibility for $u_4$ is $u_4 = (1,n-1)$ and there is no more ray, \ie $k=4$. This gives two possible fans we show on figure \ref{fig_fan0a}. 
We see these fans are congruent to the fans of figure \ref{fig_fan0b}, so $(1)$ is true.  
  
\begin{figure}[h!]
	\begin{subfigure}[t]{0.8\textwidth}
	\centering
	\begin{tikzpicture}[scale=1] 
	
	\foreach \i in {-1,...,1} {
		\foreach \j in {-2,...,2} {
			\draw node at (\i,\j) {\huge $\cdot$} ;
		}
	}
	
	\draw[gray, dashed] (0,-2) to (0,2)	;
	\draw[gray, dashed] (-1,0) to (1,0)	;
	
	\draw[->] (0,0) to (1,2) ;
	\draw[->] (0,0) to (1,1) ;
	\draw[->] (0,0) to (-1,-1) ;
	\draw[->] (0,0) to (-1,-2) ;
	
	\draw node[right] at (1,2) {\scriptsize $(1,n)$} ;
	\draw node[right] at (1,1) {\scriptsize $(1,n-1)$} ;
	\draw node[left] at (-1,-1) {\scriptsize $(-1,1-n)$} ;
	\draw node[left] at (-1,-2) {\scriptsize $(-1,-n)$} ;
	
	
	\foreach \i in {7,...,9} {
		\foreach \j in {-2,...,3} {
			\draw node at (\i,\j) {\huge $\cdot$} ;
		}
	}
	
	\draw[gray, dashed] (8,-2) to (8,3)	;
	\draw[gray, dashed] (7,0) to (9,0)	;
	
	\draw[->] (8,0) to (9,3) ;
	\draw[->] (8,0) to (9,2) ;
	\draw[->] (8,0) to (9,1) ;
	\draw[->] (8,0) to (7,-2) ;
	
	\draw node[right] at (9,3) {\scriptsize $(1,n)$} ;
	\draw node[right] at (9,2) {\scriptsize $(1,n-1)$} ;
	\draw node[right] at (9,1) {\scriptsize $(1,n-2)$} ;
	\draw node[left] at (7,-2) {\scriptsize $(-1,1-n)$} ;
	
	\end{tikzpicture}
	\caption{The two possible fans without vertical ray...}
	\label{fig_fan0a}
\end{subfigure}
\begin{subfigure}[t]{0.8\textwidth}
	\centering
	\begin{tikzpicture}[scale=1] 
	
	\foreach \i in {-1,...,1} {
		\foreach \j in {-1,...,1} {
			\draw node at (\i,\j) {\huge $\cdot$} ;
		}
	}
	
	\draw[gray, dashed] (0,-1) to (0,1)	;
	\draw[gray, dashed] (-1,0) to (1,0)	;
	
	\draw[->] (0,0) to (1,0) ;
	\draw[->] (0,0) to (0,1) ;
	\draw[->] (0,0) to (-1,0) ;
	\draw[->] (0,0) to (0,-1) ;
	
	\draw node[right] at (1,0) {\scriptsize $(1,0)$} ;
	\draw node[above] at (0,1) {\scriptsize $(0,1)$} ;
	\draw node[left] at (-1,0) {\scriptsize $(-1,0)$} ;
	\draw node[below] at (0,-1) {\scriptsize $(0,-1)$} ;
	
	
	\foreach \i in {7,...,9} {
		\foreach \j in {-1,...,1} {
			\draw node at (\i,\j) {\huge $\cdot$} ;
		}
	}
	
	\draw[gray, dashed] (8,-1) to (8,1)	;
	\draw[gray, dashed] (7,0) to (9,0)	;
	
	\draw[->] (8,0) to (7,1) ;
	\draw[->] (8,0) to (8,1) ;
	\draw[->] (8,0) to (9,1) ;
	\draw[->] (8,0) to (8,-1) ;
	
	\draw node[left] at (7,1) {\scriptsize $(-1,1)$} ;
	\draw node[above] at (8,1) {\scriptsize $(0,1)$} ;
	\draw node[right] at (9,1) {\scriptsize $(1,1)$} ;
	\draw node[below] at (8,-1) {\scriptsize $(0,-1)$} ;
	
	\end{tikzpicture}
	\caption{...are congruent to these.}
	\label{fig_fan0b}
\end{subfigure}
	\caption{}
	\label{fig_fan0}
\end{figure}
 
Assume now that $\F$ has a single vertical ray. By symmetry we can assume that this ray is generated by $(0,-1)$. Let $n\in\N$ be the biggest integer such that $\F$ has a ray generated by $(\pm 1, n)$. By symmetry we can only consider the case of $(1,n)$. Necessarily $n>0$ otherwise $\F$ is not complete. 

Let $u_1,\dots,u_k$ be the primitive generators of the rays of $\F$, taken in the anticlockwise direction and with $u_1 = (1,n)$. As previously one has $u_2 = (-1,1-n)$. The sequence of generators afterwards is of the form $u_j = (-1,3-n-j)$ for $3 \leq j \leq \ell -1$ for some $3 \leq \ell \leq k$, and $u_\ell = (0,-1)$. The $k-\ell$ remaining generators are necessarily the vectors $u_j = (1,n-k+j-1)$ for $\ell +1 \leq j \leq k$. Hence $\F$ is the fan of figure \ref{fig_fan1a}. Consider the non-singular and $h$-transverse fan $\F'$ depicted in figure \ref{fig_fan1b}. Then $\F$ can be obtained from $\F'$ performing blow-ups, and thus $(2)$ is true.
  
\begin{figure}[h!]
	\begin{subfigure}[t]{0.49\textwidth}
	\centering
	\begin{tikzpicture}[scale=1] 
	
	\foreach \i in {-1,...,1} {
		\foreach \j in {-3,...,2} {
			\draw node at (\i,\j) {\huge $\cdot$} ;
		}
	}
	
	\draw[gray, dashed] (0,-3) to (0,2)	;
	\draw[gray, dashed] (-1,0) to (1,0)	;
	
	\draw[->] (0,0) to (1,2) ;
	\draw[->] (0,0) to (1,1) ;
	\draw[->] (0,0) to (1,0) ;
	\draw[->] (0,0) to (1,-1) ;
	\draw[->] (0,0) to (1,-2) ;
	\draw[->] (0,0) to (0,-1) ;
	\draw[->] (0,0) to (-1,-3) ;
	\draw[->] (0,0) to (-1,-2) ;
	\draw[->] (0,0) to (-1,-1) ;
	
	\draw node[above] at (1,2) {\scriptsize $(1,n)$} ;
	\draw node[below] at (0,-1) {\scriptsize $(0,-1)$} ;
	\draw node[left] at (-1,-1) {\scriptsize $(-1,1-n)$} ;
	
	\end{tikzpicture}
	\caption{The fan $\F$.}
	\label{fig_fan1a}
\end{subfigure}
\begin{subfigure}[t]{0.49\textwidth}
	\centering
	\begin{tikzpicture}[scale=1] 
	
	\foreach \i in {-1,...,1} {
		\foreach \j in {-3,...,2} {
			\draw node at (\i,\j) {\huge $\cdot$} ;
		}
	}
	
	\draw[gray, dashed] (0,-3) to (0,2)	;
	\draw[gray, dashed] (-1,0) to (1,0)	;
	
	\draw[->] (0,0) to (1,2) ;
	\draw[->] (0,0) to (0,-1) ;
	\draw[->] (0,0) to (-1,-1) ;
	
	\draw node[above] at (1,2) {\scriptsize $(1,n)$} ;
	\draw node[below] at (0,-1) {\scriptsize $(0,-1)$} ;
	\draw node[left] at (-1,-1) {\scriptsize $(-1,1-n)$} ;
	
	\end{tikzpicture}
	\caption{The fan $\F'$.}
	\label{fig_fan1b}
\end{subfigure}
	\caption{}
	\label{fig_fan1}
\end{figure}
\end{demo}

If $\Delta$ and $\Delta'$ are congruent then $\GDelta(s) = G_{\Delta'}(s)$. Hence by the previous lemma we can only consider polygons having one or two horizontal sides, \ie fans having one or two vertical rays.

\begin{theo}\label{theo_nonsing}
Let $i\in\N$ and $\F$ be a non-singular and $h$-transverse fan. Let $\Delta \in D(\F)$ and $s\in\{0,\dots, \smax(\Delta)\}$. 
If 
\begin{enumerate}
	\item the fan $\F$ has two vertical rays and 
	\[ \left\{ \begin{array}{rcl}
	\Delta &>& 2(i+2) \\
	e^{+\infty}(\Delta) &>& i+ (a(\Delta)-\midaD+1)d_\F \\
	e^{-\infty}(\Delta) &>& \max(i+2s, i+ (a(\Delta)-\midaD+1)d_\F)
	\end{array} \right. \]
\end{enumerate}
or
\begin{enumerate}[resume]
	\item the fan $\F$ has a single vertical ray generated by $(0,\epsilon)$, with $\epsilon \in \{-1,1\}$, and 
	\[ \left\{ \begin{array}{rcl}
	\Delta &>& 2(i+2) \\
	a(\Delta) &>& \max(i+2s, 5(i+1)+6) \\
	e^{\epsilon\infty}(\Delta) &>& \max(i+2s, 5(i+1)+6, i+ (a(\Delta)-\midaD+1)d_\F)
	\end{array} \right. \]
\end{enumerate}
then one has
\[ \< \GDelta(s)\>_i = P_i(y(\Delta),\chi(\Delta),s) .\]
\end{theo}

\begin{demo}
In case $(1)$ we can apply theorem \ref{theo_nonsinghori}.

In case $(2)$, by lemma \ref{lemme_fan} and up to symmetries the fan $\F$ is obtained by several blow-ups from the fan of figure \ref{fig_fan1b}. 
 This fan is congruent to the fan whose rays generators are $(-1,0)$, $(0,-1)$ and $(1,1)$, which defines $\CP^2$. Hence we can use corollary \ref{coro_CP2}, invariance under $\GL_2(\Z)$ and translations, and several applications of corollary \ref{coro_blowup} to conclude.
\end{demo}

\subsection{The case of singular surfaces}\label{subsec_maintheosing}

The results of sections \ref{subsec_maintheopart} and \ref{subsec_blowupCP2} can be extended to singular surfaces. Recall that for a polygon $\Delta$, we denote by $n_k(\Delta)$ its number of vertices of index $k$.
Recall also that $Q_i$ is the degree $i$ coefficient of $ A_0^s A_1^{y-2-2s} \prod_k A_2(x^k)^{n_k}$.
Theorem \ref{theo_nonsinghori} generalizes as follows.

\begin{theo}  \label{theo_singhori} 
Let $i\in\N$. and $\F$ be a $h$-transverse fan having two vertical rays. Let $\Delta \in D(\F)$ and $s\in \{0,\dots, \smax(\Delta) \}$. If $(\Delta,s)$ satisfies 
\[ \left\{ \begin{array}{rcl}
 \Delta &>& 2(i+2) \\
 e^{-\infty}(\Delta) &>& \max(i+ (a(\Delta)-\midaD+1)d_\F,i+2s) \\
 e^{+\infty}(\Delta) &>& i+ (a(\Delta)-\midaD+1)d_\F
\end{array} \right. \]
 then 
\[ \< \GDelta(s)\>_i = Q_i(y(\Delta),s,n_1(\Delta),\dots). \]
\end{theo}

\begin{rk}
If $n_k(\Delta)=0$ unless $k=1$, the product over $k$ is just $A_2^{\chi(\Delta)}$ and we recover theorem \ref{theo_nonsinghori}.
\end{rk}

\begin{rk}
For $k\geq 2$, $n_k(\Delta)$ is the number of singularities of $X_\F$ of index $k$, see \cite[proposition A.1]{liu_severi_2018}. Hence the term taking into account the singularities $\prod_{k\geq2} A_2(x^k)^{n_k}$ is the same as the one in \cite[corollary 1.10]{liu_severi_2018}.
\end{rk}

\begin{demo} 
	For all $k\geq 2$ let $m_k=n_k(\Delta)$ and $m_1 = n_1(\Delta)-4$. 	For a vertex of index $k \geq 1$ of $\Delta$, the adjacent edges have direction vectors $(q,-1)$ and $(q\pm k,-1)$ for some $q\in\Z$. Hence in the choice of the function $r$ for the construction of a diagram $\D$  (see proof of lemma \ref{lemme_rgamma}) the vector $\tilde\gamma$ will contribute to the codegree  by
	\[ \dsum_{j=1}^i k(\tilde\gamma_j -j) \]
	and  if we want $\codeg(\D) \leq i$ then the corresponding vector $\gamma$ will satisfy 
	\[ \codeg(\gamma) \leq \dfrac{i}{k} . \]
	Thereby, to choose the functions $r$ and $\ell$ we need, for every $k\geq 1$, to choose $m_k$ vectors $\gamma_k^1,\dots, \gamma_k^{m_k}$ of codegree at most $i/k$ and such that 
	\[  \dsum_{k\geq 1} \dsum_{j=1}^{m_k} \codeg(\gamma_k^j) \leq i . \]
	If the divergence contributes $c$ to the codegree of $\D$, then by lemma \ref{lemme_cardBi} the number of choices for $(r,\ell)$ is
	\[ \dsum_{c_1 + c_2 + \dots = c}\ \ \dprod_{k \geq 1} \ \ \dsum_{j_1+\dots + j_{m_k} = c_k/k} p(j_1) \dots p(j_{m_k}) \]
	where $c_k$ is the contribution of the vertices of index $k$.
	The corresponding generating series is
	\[  \dprod_{k^\geq1} A_2(x^k)^{m_k} \]
	and the rest of the proof is as in theorem \ref{theo_nonsinghori}.
\end{demo}	

If $\F$ is a fan having two consecutive rays primitively generated by $u$ and $v$, we construct a new fan $\tilde\F$ from $\F$ by adding a ray generated by $mu+v$ for some $m \geq 1$. In other words, we replace a $2$-dimensional cone of $\F$ of index $|\det(u,v)|$ by two $2$-dimensional cones of indices $|\det(u,v)|$ and $m$.
Proposition \ref{prop_blowup} deals with the case $m=1$ and generalizes as follows. Recall the conditions 
\[ \stari \left\{ \begin{array}{rcl}
 \Delta &>& 2(i+2) \\
 e^{-\infty}(\Delta) &>& i+2s \\
 e^{-\infty}(\Delta) &>& i+ (a(\Delta)-\midaD+1)d_\F
\end{array} \right. . \] 

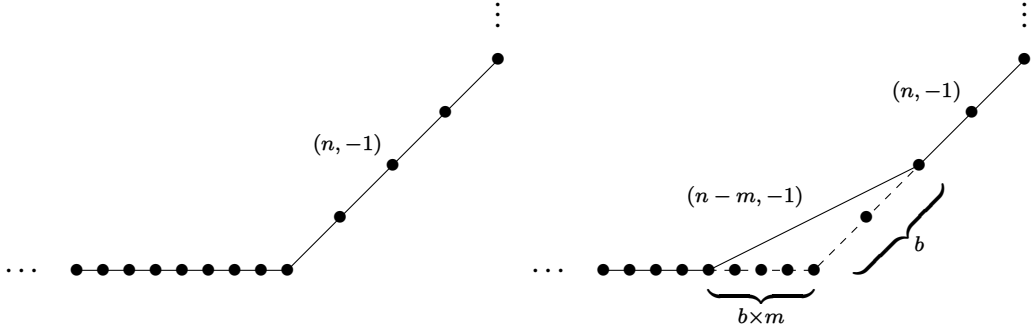
\begin{figure}[h] 
		\centering
	\begin{tikzpicture}[scale=0.7]
		
		\draw (-4,0) to (0,0) to node[above left] {\scriptsize $(n,-1)$} (4,4) ;
		\foreach \p in {(-4,0),(-3.5,0),(-3,0),(-2.5,0),(-2,0),(-1.5,0),(-1,0),(-0.5,0), (0,0), (1,1),(2,2),(3,3),(4,4)}
			{\draw node at \p {$\bullet$} ;} 
		\draw node at (-5,0) {$\dots$} ;
		\draw node at (4,5) {$\vdots$} ;

		\draw (6,0) to (8,0) to node[above left] {\scriptsize $(n-m,-1)$} (12,2) to node[above left] {\scriptsize $(n,-1)$} (14,4);
		\draw[dashed] (8,0) to (10,0) to (12,2) ;
		\foreach \p in {(6,0),(7,0),(8,0),(9,0),(6.5,0),(7.5,0),(8.5,0),(9.5,0), (10,0),(11,1),(12,2),(13,3),(14,4)}
			{\draw node at \p {$\bullet$} ;} 
		\draw node at (5,0) {$\dots$} ;
		\draw node at (14,5) {$\vdots$} ;
		
		\draw node at (9,-0.7) {$\underbrace{\ \ \ \ \ \ \ \ \ \ \ }_{b\times m}$} ;

		\draw (11,1.4) node[right,rotate=-45] {\scriptsize $ \left. \begin{array}{c}
 \\ \\ \\  \\ \\ 		
\end{array}		 \right\}$} ;	
		\draw node at (12,0.5) {\scriptsize $b$} ;

	\end{tikzpicture}
	\caption{The polygons $\Delta$ (left) and $\tilde\Delta$ (right).} 
	\label{fig_blowupsing} 
\end{figure}

\begin{prop} \label{prop_blowupsing}
 Let $\Delta$ and $\tilde\Delta$ be the polygons whose bottom right corners are depicted in figure \ref{fig_blowupsing}. 
 There exists a series $B$ such that if $(\Delta,s)$ satisfies $\stari$ then $\< \GDelta(s) \>_i$ is given by the degree $i$ coefficient of 
 \[ A_0^s A_1^{e^{-\infty}(\Delta)-2s} A_2^2 \times B ,\]
 and if $(\tilde\Delta,s)$ satisfies $\stari$ then $\< G_{\tilde\Delta}(0;s) \>_i$ is given by the degree $i$ coefficient of 
 \[ A_0^s A_1^{e^{-\infty}(\tilde\Delta)-2s} A_2^2 A_2(x^m) \times B .\] 
\end{prop}
  
\begin{demo}
It is the same proof as proposition \ref{prop_blowup}, except that to pass from a floor diagram with Newton polygon $\Delta$ to a floor diagram with Newton polygon $\tilde\Delta$ we drop the divergence of some vertices by $m$. Hence, in the computations for $\tilde\Delta$ we need to consider $\gamma \in B_{k/m}$, where $B_{k/m} = \vide$ if $m$ does not divide $k$. We will end up with
\begin{align*}
 \< G_{\tilde\Delta}(0;s) \>_i
  &= \dsum_{k=0}^i p\left(\dfrac{k}{m}\right) \dsum_{j=0}^{i-k} \alpha_j(e^{-\infty}-bm,s) \beta_{i-k-j}(a_n)
\end{align*}
which shows that $\< G_{\tilde\Delta}(0;s) \>_i$ is the degree $i$ coefficient of the product of
\begin{itemize}
	\item the generating series of $(p(k/m))_k$ which is $A_2(x^p)$ by definition, 
	
	\item the generating series of $(\alpha_k(e^{-\infty}-bm,s))_k$ which is $A_0^s A_1^{e^{-\infty}-bm-2s} A_2^2$ by lemma \ref{lemme_binom},
	
	\item the generating series of $(\beta_k(a_n))_k$ which is $B$.
\end{itemize}
Since $e^{-\infty}(\tilde\Delta) = e^{-\infty}-bm$ we conclude.
\end{demo}

We could try to generalize theorem \ref{theo_nonsing} to singular surfaces, but some difficulties appear.

\begin{itemize}
\item If $\Delta$ is a $h$-transverse polygon with one horizontal edge, it can be obtained by several transformations described in proposition \ref{prop_blowupsing} starting from a polygon congruent to the triangle defining a weighted projective plane $\CP(1,1,n)$. Hence we can link the generating series for $\Delta$ with the one for $\CP(1,1,n)$, but still we do not know how to compute explicitly the latter (although we conjecture it fits in the general form of theorem \ref{theo_singhori}).

\item For a $h$-transverse polygon without horizontal edge, there is no reason for it to be congruent to a $h$-transverse polygon with one horizontal edge (leading to the previous case) or with two horizontal edges (leading to the case of theorem \ref{theo_singhori}). Hence we are not able to say anything for these polygons, but we conjecture they fit in the general form of theorem \ref{theo_singhori}.
\end{itemize}

\section{Removing the denominators} \label{sec_removedenom}

\subsection{Motivations}

In what precedes, computations are made and formulas are shown for tropical refined invariants as defined in theorem \ref{def_ITR}. In particular, the multiplicities we considered are products of terms of the form
\[ [n](q) = \dfrac{q^{n/2} - q^{-n/2}}{q^{1/2} - q^{-1/2}} . \]
 However some results suggest that the interesting part of the multiplicities is the numerators. For instance, \cite[theorem 5.9]{mikhalkin_quantum_2017} and \cite[theorem 4.12]{blomme_refined_2023} establish a link between tropical refined invariants without denominators and  a refined count of real curves ; \cite[theorem 5]{bousseau_tropical_2019} relates tropical refined invariants without denominators to a generating series of log Gromov-Witten invariants.  Therefore, we briefly discuss in this section what happens for our results if we get rid of the denominators in the multiplicities.

\subsection{The result without denominators.}

Let $(\D,m)$ be a marked floor diagram and $S$ be a pairing compatible with $m$. Recall the notations $E_0$, $E_1$ and $E_2$ used in definition \ref{def_mult}. We consider the partitions 
\[ E_0 = E_0^0 \sqcup E_0^\infty,\ E_1 = E_1^0 \sqcup E_1^\infty \text{ and } E_2 = E_0^{00} \sqcup E_2^{0\infty} \sqcup E_2^{\infty \infty}, \]
 where the superscript $0$ stands for bounded and the superscript $\infty$ stands for infinite, for instance $E_2^{0 \infty} = \{\{e,e'\} \in E_2\ |\ e \text{ bounded and }\ e' \text{ infinite}\}$.
  For an integer $n$ and a formal variable $t$ we set the notation $t_n = 1-t^n$. We introduce the following multiplicity :
\begin{multline*}  \mu_S^\star(\D,m)(t) = t^{\codeg(\D)} 
   \dprod_{e\in E_0^0} t_{{w(e)}}^2
   \dprod_{e\in E_0^\infty} t_1
   \dprod_{e\in E_1^0} t_{2w(e)}
   \dprod_{e\in E_1^\infty} (1+t) \\
   \dprod_{\{e,e'\}\in E_2^{00}} t_{w(e)} t_{w(e')} t_{w(e)+w(e')}
   \dprod_{\{e,e'\} \in E_2^{0 \infty}} t_{w(e)} t_{w(e)+1}
   \dprod_{\{e,e'\} \in E_2^{\infty \infty}} t_2 
\end{multline*}
and consider 
\[ \GDelta^\star(S)(t) = \dsum_{(\D,m)} \mu_S^\star(\D,m)(t) \]
where the sum runs over the isomorphism classes of marked floor diagrams with Newton polygon $\Delta$.

For $A(q) \in \C[q^{1/2},q^{-1/2}]$ a symmetric Laurent polynomial we set 
\[ \tilde A(t) = t^{\deg(A)} A(t) \in \C[t] \]
 \ie we dispel the negative powers. If $A(q) = a_0q^d + a_1 q^{d-1} + \dots + a_0 q^{-d}$, then $\tilde A(t) = a_0 + a_1 t + a_2 t^2 + \dots$ : the coefficient $a_i$ of codegree $i$ in $A(q)$ becomes the coefficient of degree $i$ of $\tilde A(t)$. Note also that $\tilde{AB} = \tilde A \tilde B$ for $A,B \in \C[q^{1/2},q^{-1/2}]$.

\begin{lemme} \label{lemme_newmult}
Let $\Delta$ be a $h$-transverse polygon. Let $s\in\{0,\dots,\smax(\Delta)\}$ and $S$ be a pairing of order $s$. One has
\[ \tilde\GDelta(s)(t) = A_0^s A_1^{y(\Delta)-2-2s} \GDelta^\star(S)(t) . \]
In particular $\GDelta^\star(S)(t)$ does not depend on the choice of the pairing $S$ of order $s$.
\end{lemme}

\begin{demo}
For any diagram $\D$ one has
\[\deg(\GDelta(s)) = \deg(\D) + \codeg(\D) = \deg(\mu_S(\D,m)) + \codeg(\D), \] 
hence 
\[ \tilde\GDelta(s)(t) = \dsum_{(\D,m)} t^{\codeg(\D)} \tilde\mu_S(\D,m) . \]
For quantum integers one has
\[
 \tilde{[n]}(t) = \dfrac{1-t^n}{1-t}, \ 
 \tilde{[n]^2}(t) = \left(\dfrac{1-t^n}{1-t}\right)^2, \text{ and }
 \tilde{[n]_2}(t) = \dfrac{1-t^{2n}}{1-t^2}
\]
and we get
\[  \tilde{\mu_S}(\D,m)(t) =
	\dprod_{e \in E_0} \left(\dfrac{1-t^{w(e)}}{1-t} \right)^2 
	\dprod_{e\in E_1} \dfrac{1-t^{2w(e)}}{1-t^2} 
	\dprod_{\{e,e'\} \in E_2} \dfrac{(1-t^{w(e)}) (1-t^{w(e')}) (1-t^{w(e)+w(e')})}{(1-t)^2 (1-t^2)}   .\]
The product over $E_0$ is rewritten 
\[  \dprod_{e\in E_0^0} \left(\dfrac{1-t^{w(e)}}{1-t}\right)^2 \dprod_{e\in E_0^\infty} \dfrac{1-t}{1-t} 
=
 A_1^{2|E_0^0| + |E_0^\infty|} \dprod_{e\in E_0^0} t_{{w(e)}}^2 \dprod_{e\in E_0^\infty} t_1, \]
the one over $E_1$ becomes
\[ \dprod_{e\in E_1^0} \dfrac{1-t^{2w(e)}}{1-t^2} \dprod_{e\in E_1^\infty} \dfrac{(1-t)(1+t)}{1-t^2} 
=
 A_0^{|E_1|} A_1^{-|E_1^\infty|} \dprod_{e\in E_1^0} t_{2w(e)} \dprod_{e\in E_1^\infty} (1+t), \]
and the one over $E_2$ gives
\begin{multline*}
	\dprod_{\{e,e'\}\in E_2^{00}}  \dfrac{(1-t^{w(e)})(1-t^{w(e')})(1-t^{w(e)+w(e')})}{(1-t)^2 (1-t^2)} 
	\dprod_{\{e,e'\}\in E_2^{0 \infty}} \dfrac{(1-t^{w(e)})(1-t^{w(e)+1})}{(1-t)(1-t^2)} 
	\dprod_{\{e,e'\}\in E_2^{\infty \infty}}\dfrac{1-t^2}{1-t^2} \\
= A_0^{|E_2|} A_1^{2|E_2^{00}| + |E_2^{0\infty}|}
	\dprod_{\{e,e'\}\in E_2^{00}} t_{w(e)} t_{w(e')} t_{w(e)+w(e')}
	\dprod_{\{e,e'\} \in E_2^{0 \infty}} t_{w(e)} t_{w(e)+1}
	\dprod_{\{e,e'\} \in E_2^{\infty \infty}} t_2 .
\end{multline*}
Putting together the three products, the total power for $A_0$ is $|E_1|+|E_2| =s$. Given the system
\[ \left\{ \begin{array}{rcl}
 |E_1|+|E_2| &=& s  \\
 |V|-|E^b| &=& 1 \\
 2|V|+|E^\infty| &=& y(\Delta) \\
 |E^\alpha| &=& |E_0^\alpha| + |E_1^\alpha| + 2|E_2^{\alpha \alpha}| + |E_2^{\alpha \beta}| , \text{ where } \{\alpha,\beta\} = \{0,\infty\}
\end{array}  \right. \]
the total power for $A_1$ is $y(\Delta)-2-2s$. Hence we eventually get 
\[ \tilde\GDelta(s)(t) = A_0^s A_1^{y(\Delta)-2-2s} \dsum_{(\D,m)} \mu_S^\star(\D,m)(t) = A_0^s A_1^{y(\Delta)-2-2s} \GDelta^\star(S)(t) .\]
\end{demo}

\begin{theo}
Let $\F$ be a $h$-transverse fan and $\Delta \in D(\F)$. Let $s\in\{0,\dots,\smax(\Delta)\}$ and $S$ be a pairing of order $s$. 
\begin{enumerate}
	\item If $\F$ is non-singular, let $i_m$ be the maximal integer such that $(\Delta,s)$ satisfies one of the conditions $(1)$ or $(2)$ of theorem \ref{theo_nonsing}. One has
\[ \GDelta^\star(S)(t) = A_2^{\chi(\Delta)} \mod t^{i_m} . \]
In particular, it does not depend on $S$ neither on the order $s$ of $S$.
	
	\item If $\F$ has two vertical rays, let $i_m$ be the maximal integer such that $(\Delta,s)$ satisfies the conditions of theorem \ref{theo_singhori}. One has
\[ \GDelta^\star(S)(t) = \dprod_{k\geq1} A_2(x^k)^{n_k(\Delta)} \mod t^{i_m} .\]
In particular, it does not depend on $S$ neither on the order $s$ of $S$.
\end{enumerate}
\end{theo}

\begin{demo}
\begin{enumerate}
	\item Theorem \ref{theo_nonsing} shows that 
\[ \tilde\GDelta(s)(t) = A_0^s A_1^{y(\Delta)-2-2s} A_2^{\chi(\Delta)} \mod t^{i_m}  \]
\ie  
\[ A_0^s A_1^{y(\Delta)-2-2s} \GDelta^\star(S)(t) = A_0^s A_1^{y(\Delta)-2-2s} A_2^{\chi(\Delta)} \mod t^{i_m}  \]
hence
\[ \GDelta^\star(S)(t) = A_2^{\chi(\Delta)} \mod t^{i_m} . \]

	\item Theorem \ref{theo_singhori} shows that
\[ \tilde\GDelta(s)(t) = A_0^s A_1^{y(\Delta)-2-2s} \dprod_{k\geq1} A_2(x^k)^{n_k(\Delta)} \mod t^{i_m} \]
which leads to
\[ \GDelta^\star(S)(t) = \dprod_{k\geq1} A_2(x^k)^{n_k(\Delta)} \mod t^{i_m} . \]
\end{enumerate}
\end{demo}

In other words, in case $(1)$ the degree $i$ coefficient of $\GDelta^\star(s)$ is given by a polynomial $P_i^\star \in \N[\chi]$ of degree $i$, and the generating series of $(P_i^\star)_i$ is $A_2^\chi$.
In case $(2)$, the degree $i$ coefficient of $\GDelta^\star(s)$ is given by a polynomial $Q_i^\star \in \N[n_1,n_2,\dots]$ and the generating series of $(Q_i^\star)_i$ is $\prod_{k\geq1} A_2(x^k)^{n_k}$.
These expressions for the generating series are easier than the ones of theorems \ref{theo_nonsing} and \ref{theo_singhori}. In particular, it is remarkable that they do not depend neither on $y$ nor on $s$.

\section{Some technical lemmas}\label{sec_lemmas}

In this section we prove few lemmas regarding the quantities we handled in section \ref{sec_univseries}. We refer to the beginning of section \ref{subsec_maintheopart} for the notations.
 
\begin{lemme} \label{lemme_cardBi}
Let $i\in\N$. The cardinality of $B_i$ is $p(i)$, the number of partitions of $i$.
\end{lemme}

\begin{demo}
Let $\lambda = (\lambda_1,\dots,\lambda_n)$ be a partition of $i$. For $1\leq k\leq i$ let $w_k(\lambda)$ be the number of $k$ in the sequence $\lambda$. Then $(w_1(\lambda),\dots,w_i(\lambda)) \in B_i$. Conversely, any $w\in B_i$ defines a partition $\lambda(w) = (\lambda_1,\dots,\lambda_n)$ where $\lambda_j=k$ for $w_i+\dots + w_{k+1}+1 \leq j \leq w_i + \dots + w_k$.
\end{demo}

\begin{lemme} \label{lemme_binom}
\begin{enumerate}
	\item Let $a,p\in\N$ and $S\in\N^\N$. One has 
	\[\calN(a,p,S) =  A_1^a A_2^p \dprod_{k\geq0} (x^{2k})^{S_k}.\]

	\item Let $s,a,p\in\N$. One has 
	\[ \dsum_{S\in d(s)} \binom{s}{S} \calN(a,p,S) = A_0^s A_1^a A_2^p. \]
\end{enumerate}
\end{lemme}

\begin{demo}
\begin{enumerate}
	\item We will prove the following formula by induction over $n$ :
\[ 
  \calN(a,p,S) =  A_1^a  \times \dprod_{k=1}^{n-1} \dfrac{(x^k)^{2S_k}}{(1-x^k)^p} \times (1-x^n)^{a+(n-1)p}
  \times \calN_n(a,p,S)
	\]
where 
\[  
\calN_n(a,p,S) = \dsum_{(u_n,\dots)} \nu_{\geq n}(a,p,u,S) x^{\codeg_n(u)} \left(\dfrac{1-x}{1-x^n}\right)^{\som_n(u-2S)} .
\]	
Note that $\calN_1=\calN$ so the formula is true for $n=1$. Assume that the formula holds for some $n\geq1$. Writing $\nu_\star$ as a shortcut for $\nu_\star(a,p,u,S)$ and using the binomial formula we get
\begin{align*}
 & \calN_n(a,p,S) \\ 
 &= \dsum_{(u_{n+1},\dots)} \nu_{\geq n+1} x^{\codeg_{n+1}(u)} \left(\dfrac{1-x}{1-x^{n}}\right)^{\som_{n+1}(u-2S)}  \dsum_{u_{n}\geq0}  \nu_{n} \left( \dfrac{x^{n} - x^{n+1}}{1-x^{n}} \right)^{u_n-2S_n} (x^n)^{2S_n} \\ 
    &= \dsum_{(u_{n+1},\dots)} \nu_{\geq n+1} x^{\codeg_{n+1}(u)} \left(\dfrac{1-x}{1-x^{n}}\right)^{\som_{n+1}(u-2S)}  \left( \dfrac{1-x^{n+1}}{1-x^{n}} \right)^{a+np-\som_{n+1}(u-2S)} (x^n)^{2S_n} \\   
    &= (x^n)^{2S_n} \left( \dfrac{1-x^{n+1}}{1-x^n} \right)^{a+np} \calN_{n+1}(a,p,S) 
\end{align*}
which gives the desired formula for $\calN(a,p,S)$.

	\item This results from the previous point and the multinomial formula :
	\[   \dsum_{S\in d(s)} \binom{s}{S} \calN(a,p,S) =  A_1^a A_2^p \dsum_{S\in d(s)} \binom{s}{S} \dprod_{k\geq0} (x^{2k})^{S_k} = A_0^s A_1^a A_2^p . \]
\end{enumerate}
\end{demo}

It is shown in \cite[lemma 3.5]{brugalle_polynomiality_2022} that $\Phi_\ell$ is a polynomial of degree $\ell$ in $k$. Actually we can give an explicit formula for this polynomial.

\begin{lemme}\label{lemme_Phi}
For $k,\ell \in\N$ one has 
\[ \Phi_\ell(k) = \binom{2k+\ell-1}{\ell} . \]
\end{lemme}

\begin{demo}
Here is a combinatorial proof. By definition, 
\[ \Phi_\ell(k) = \dsum_{\substack{i_1 + \dots + i_k = k+\ell \\ i_j \geq 1 }} \dprod_{j=1}^k i_j = \dsum_{\substack{i_1 + \dots + i_k = \ell \\ i_j \geq 0 }} \dprod_{j=1}^k (i_j+1).  \]
The interpretation is the following. We first need to choose a decomposition of $\ell$ into $k$ parts. This is the same as taking $k+\ell-1$ white boxes in line and choosing $k-1$ of them to be black. The black boxes determine $k$ groups of $i_1,\dots,i_k$ white boxes. We add one blue box in each of these groups. Then $\prod_{j=1}^k(i_j+1)$ is the number of ways of choosing the places of the blue boxes. 

This is the same as taking $2k+\ell-1$ whites boxes, then choosing and coloring $2k-1$ of them alternately blue and black. Thus 
\[\Phi_\ell(k) = \dbinom{2k+\ell-1}{2k-1} = \dbinom{2k+\ell-1}{\ell}.\]
\end{demo}

For $\ell \in \N$ we can thus extend $\Phi_\ell$ to a function $t \in \R \mapsto \Phi_\ell(t)$.

\begin{coro} \label{coro_genPhi}
Let $t\in\R$. The generating series of $(\Phi_\ell(t))_{\ell\in\N}$ is $A_1^{2t}$. 
\end{coro}

\begin{demo}
This comes from the fact that $\binom{2t+\ell-1}{\ell} = (-1)^\ell \binom{-2t}{\ell}$ (see the upper negation in \cite[table 174]{graham_concrete_1994}) and from the binomial formula.
\end{demo}

\begin{lemme}[{\cite[corollary 3.6]{brugalle_polynomiality_2022}}] \label{lemme_mult}
Let $i,k \in \N$ and $a_1,\dots,a_k >i$ be integers. One has
\[  \left\< \dprod_{j=1}^k [a_j]^2 \right\>_i = \Phi_i(k).  \]
\end{lemme}


\renewcommand{\refname}{References}
\bibliographystyle{alpha}
\bibliography{ref}

\begin{thebibliography}{GKP94}

\bibitem[AB13]{ardila_universal_2013}
Federico Ardila and Florian Block.
\newblock Universal polynomials for {S}everi degrees of toric surfaces.
\newblock {\em Advances in Mathematics}, 237:165--193, 2013.

\bibitem[BCK14]{block_computing_2014}
Florian Block, Susan~Jane Colley, and Gary Kennedy.
\newblock Computing {Severi} {Degrees} with {Long}-edge {Graphs}.
\newblock {\em Bulletin of the Brazilian Mathematical Society, New Series},
  45(4):625--647, 2014.

\bibitem[BG16]{block_refined_2016}
Florian Block and Lothar Göttsche.
\newblock Refined curve counting with tropical geometry.
\newblock {\em Compositio Mathematica}, 152(1):115--151, January 2016.

\bibitem[BJP22]{brugalle_polynomiality_2022}
Erwan Brugallé and Andrés Jaramillo~Puentes.
\newblock Polynomiality properties of tropical refined invariants.
\newblock {\em Combinatorial Theory}, 2, 2022.

\bibitem[Blo11]{block_computing_2011}
Florian Block.
\newblock Computing node polynomials for plane curves.
\newblock {\em Mathematical Research Letters}, 18(4):621--643, 2011.

\bibitem[Blo23]{blomme_refined_2023}
Thomas Blomme.
\newblock Refined count of oriented real rational curves.
\newblock {\em Journal of Algebraic Geometry}, May 2023.

\bibitem[BM07]{brugalle_enumeration_2007}
Erwan Brugallé and Grigory Mikhalkin.
\newblock Enumeration of curves via floor diagrams.
\newblock {\em Comptes Rendus Mathématiques de l'Académie des Sciences de
  Paris}, 345(6):329--334, September 2007.

\bibitem[BM08]{brugalle_floor_2008}
Erwan Brugallé and Grigory Mikhalkin.
\newblock Floor decompositions of tropical curves : the planar case.
\newblock {\em Proceedings of 15th Gökova Geometry-Topology Conference}, pages
  64--90, 2008.

\bibitem[Bou19]{bousseau_tropical_2019}
Pierrick Bousseau.
\newblock Tropical refined curve counting from higher genera and lambda
  classes.
\newblock {\em Invent. math}, 215:1 -- 79, 2019.

\bibitem[DFI95]{di_francesco_quantum_1995}
Philippe Di~Francesco and Claude Itzykson.
\newblock Quantum intersection rings.
\newblock In {\em The Moduli Space of Curves}, pages 81--148, Boston, MA, 1995.
  Birkh{\"a}user Boston.

\bibitem[FM10]{fomin_labeled_2010}
Sergey Fomin and Grigory Mikhalkin.
\newblock Labeled floor diagrams for plane curves.
\newblock {\em Journal of the European Mathematical Society}, pages 1453--1496,
  2010.

\bibitem[GKP94]{graham_concrete_1994}
Ronald~Lewis Graham, Donald~Ervin Knuth, and Oren Patashnik.
\newblock {\em Concrete mathematics: a foundation for computer science}.
\newblock Addison-Wesley, Reading, Mass, 2nd ed edition, 1994.

\bibitem[GS19]{gottsche_refined_2019}
Lothar Göttsche and Franziska Schroeter.
\newblock Refined broccoli invariants.
\newblock {\em Journal of Algebraic Geometry}, 28(1):1--41, 2019.

\bibitem[Gö98]{gottsche_conjectural_1998}
Lothar Göttsche.
\newblock A conjectural generating function for numbers of curves on surfaces.
\newblock {\em Communications in Mathematical Physics}, 196(3):523--533,
  September 1998.

\bibitem[KP04]{kleiman_node_2004}
Steven~Lawrence Kleiman and Ragni Piene.
\newblock Node polynomials for families: methods and applications.
\newblock {\em Mathematische Nachrichten}, 271:69--90, 2004.

\bibitem[KST11]{kool_short_2011}
Martijn Kool, Vivek Shende, and Richard Thomas.
\newblock A short proof of the {Göttsche} conjecture.
\newblock {\em Geometry \& Topology}, 15(1):397--406, March 2011.

\bibitem[Liu16]{liu_combinatorial_2016}
Fu~Liu.
\newblock A combinatorial analysis of {Severi} degrees.
\newblock {\em Advances in Mathematics}, 298:1--50, August 2016.

\bibitem[LO18]{liu_severi_2018}
Fu~Liu and Brian Osserman.
\newblock Severi degrees on toric surfaces.
\newblock {\em Journal f{\"u}r die reine und angewandte Mathematik (Crelles
  Journal)}, (739):121--158, 2018.

\bibitem[Mik05]{mikhalkin_enumerative_2005}
Grigory Mikhalkin.
\newblock Enumerative tropical algebraic geometry in {$\R^2$}.
\newblock {\em Journal of the American Mathematical Society}, 18(2):313--377,
  January 2005.

\bibitem[Mik17]{mikhalkin_quantum_2017}
Grigory Mikhalkin.
\newblock {Quantum indices and refined enumeration of real plane curves}.
\newblock {\em Acta Mathematica}, 219(1):135 -- 180, 2017.

\bibitem[Tze12]{tzeng_proof_2012}
Yu-jong Tzeng.
\newblock A proof of the {G}{\"o}ttsche-{Y}au-{Z}aslow formula.
\newblock {\em Journal of Differential Geometry}, 90:439--472, 2012.

\bibitem[Vai95]{vainsencher_enumeration_1995}
Israel Vainsencher.
\newblock Enumeration of {$n$}-fold tangent hyperplanes to a surface.
\newblock {\em Journal of Algebraic Geometry}, 4(3):503--526, 1995.

\end{thebibliography}

\end{document}